\documentclass{amsart}
\usepackage{amsmath,amsthm,bbm}
\usepackage{amsfonts,amssymb}
\usepackage[english]{babel}

 \newtheorem{theorem}{Theorem}[section]

 \newtheorem{proposition}[theorem]{Proposition}
 \newtheorem{definition}[theorem]{Definition}
 
 \newtheorem{remark}[theorem]{Remark}
 
 \newcounter{figures}[section]

\def\bC{{\mathbb C}}

\def\NN{{\mathbb N}}
\def\PP{{\mathbb P}}
\def\RR{{\mathbb R}}
\def\SS{{\mathbb S}}
\def\ZZ{{\mathbb Z}}

  \def\cF{\mathcal{F}}
  
  \def\cH{\mathcal{H}}

  \def\cO{\mathcal{O}}

  \def\cX{\mathcal{X}}

\def\CF{{\mathcal F}}
\def\CH{{\mathcal H}}

\def\CK{{\mathcal K}}
\def\CL{{\mathcal L}}
\def\CM{{\mathcal M}}

\def\CV{{\mathcal V}}

\def\proj{\operatorname{Proj}}

\def\supp{\operatorname{supp}}

\def\eps{{\varepsilon}}
\def\a{{\alpha}}
\def\b{{\beta}}
\def\k{{\kappa}}
\def\g{{\gamma}}
\def\la{{\lambda}}
\def\ha{{\hat a}}

\def\ph{{\varphi}}

\def\ONE{{\mathbbm 1}}

\def\cc{\tilde{c}}
\def\CC{c^\diamond}
\def\Q{{\mathcal Q}}
\def\PP{{\mathsf {P}}}
\def\TT{\tilde{T}}
\def\tP{\tilde{P}}

\def\ph{{\varphi}}
\newcommand{\wh}{\hat}
\def\xxi{{\xi}}
\def\dd{{\rho}}
\def\ddT{{\rho}}
\def\zz{{z}}
\def\kk{{\kappa}}

\begin{document}

\title{Sub-exponentially localized kernels and frames induced by orthogonal expansions}

\author{Kamen Ivanov, Pencho Petrushev and Yuan Xu}

\address{Institute of Mathematics and Informatics\\
Bulgarian Academy of Sciences\\ 1113 Sofia\\ Bulgaria}
\email{ivanov@math.sc.edu}

\address{Department of Mathematics\\University of South Carolina\\
Columbia, SC 29208\\
and Institute of Mathematics and Informatics, Bulgarian Academy of Sciences}
\email{pencho@math.sc.edu}

\address{Department of Mathematics\\ University of Oregon\\
Eugene, Oregon 97403-1222.}\email{yuan@math.uoregon.edu}

\date{September 18, 2008}




\thanks{The second author has been supported by NSF Grant DMS-0709046
and the third author by NSF Grant DMS-0604056.}



\begin{abstract}
The aim of this paper is to construct sup-exponentially localized kernels and frames in the context
of classical orthogonal expansions, namely, expansions in
Jacobi polynomials, spherical harmonics,
orthogonal polynomials on the ball and simplex,
and Hermite and Laguerre functions.
\end{abstract}

\maketitle

\pagestyle{myheadings}
\thispagestyle{plain}
\markboth{KAMEN IVANOV, PENCHO PETRUSHEV AND YUAN XU}
         {LOCALIZED KERNELS AND FRAMES INDUCED BY ORTHOGONAL EXPANSIONS}

\section{Introduction}\label{introduction}
\setcounter{equation}{0}

Orthogonal expansions have been recently used for the construction of kernels and frames (needlets)
with localization that is faster than the reciprocal of any polynomial rate
in non-standard settings such as on
the sphere, interval and ball with weights, and in the context of Hermite and Laguerre expansions.
The main purpose of this article is to show that the rapid decay of that sort of
kernels and needlets can be improved to sub-exponential.
In order to best present our results, it is perhaps suitable
first to exhibit and illustrate the main principles 
and ideas which guided us in this undertaking.

\subsection{Localization principle}
\label{localizatin-principle}

Here we briefly revisit the ``localization principle" described in \cite{PX2}.
Let $(E, \mu)$ be a measure space with $E$ a metric space
and suppose that there is an orthogonal decomposition
$
L^2(E, \mu) = \bigoplus_{n=0}^\infty \CV_n,
$
where $\CV_n$ are finite dimensional subspaces.
Let $P_n$ be the kernel of the orthogonal projector
$\proj_n: L^2(E, \mu)\to\CV_n$, i.e.
$$
(\proj_n f)(x) =\int_E P_n(x,y)f(y)d\mu, \quad f\in L^2(E, \mu).
$$
We are interested in kernels of the form
\begin{equation} \label{def:Ln-gen}
L_n(x,y) := \sum_{j=0}^\infty \ha\Big(\frac{j}{n}\Big) P_j(x,y),
\end{equation}
where the cutoff function $\ha$ is compactly supported and in $C^\infty$.
For all our purposes it suffices to only consider cutoff functions obeying
the following definition:


\begin{definition}\label{cutoff-d1}
A function $\ha\in C^\infty [0, \infty)$ $($$\ha \ge 0$ if needed$)$ is said to be admissible of
type $(a)$, $(b)$ or $(c)$ if it obeys the conditions:

$(a)$ $\supp \ha \subset [0, 2]$ and $\ha(t) = 1$, $t\in [0, 1]$; or

$(b)$ $\supp \ha \subset [1/2, 2]$; or

$(c)$ $\supp \ha \subset [1/2, 2]$ and $|\ha(t)|^2+|\ha(t/2)|^2=1$ for $t\in [1, 2]$.
\end{definition}
Note that (c) is a subcase of (b). We list it separately as the
additional requirement in (c) plays an essential role in our
construction of tight frames.

The {\bf localization principle} put forward in \cite{PX2} says that
for all ``natural orthogonal systems" the kernels $\{L_n(x,y)\}$ decay
at rates faster than any inverse polynomial rate away from the main
diagonal $y=x$ in $E\times E$ with respect to the distance in $E$.
This principle is very well-known in the case of the trigonometric system
(and the Fourier transform) and not so long ago was established for
spherical harmonics \cite{NPW1, NPW2},
Jacobi polynomials \cite{BD,PX1},
orthogonal polynomials on the ball \cite{PX2}, and
Hermite and Laguerre functions \cite{Dzub, Epp2, PX3, KPPX}.

Surprisingly, however, the localization principle as formulated above
fails to be true for tensor product Jacobi polynomials and, in particular,
for tensor product Legendre or Chebyshev polynomials,
as will be shown in \S\ref{not-valid-principle}.

To grasp the notion of {\em rapidly decaying kernels} of form (\ref{def:Ln-gen})
let us illustrate them in the simple case of Chebyshev polynomials.
Denote by $\tilde{T}_n$, $n=0, 1, \dots$, the weighted-$L^2$-normalized Chebyshev polynomials
of first kind. Then for an admissible cutoff function $\ha$ the kernels from
(\ref{def:Ln-gen}) take the form
\begin{equation}\label{Chebyshev-Kernel-intr}
L_n(x,y) := \sum_{j=0}^\infty \ha\left(\frac{j}{n}\right) \tilde{T}_j(x)\tilde{T}_j(y)
\end{equation}
and satisfy (see \S\ref{Examples})
\begin{equation}\label{Loc-Chebyshev-intr}
|L_n(x,y)|\le c_\sigma n(1+n\dd(x, y))^{-\sigma}, \quad x, y \in [-1, 1],
\end{equation}
for arbitrarily large $\sigma>0$ but the constant $c_\sigma$ depends on sigma;
here $\dd$ is the distance
$\dd(x, y):=|\arccos x - \arccos y|$.
The above estimate suggests that the localization of $L_n(x,y)$ can eventually be improved to
a localization of exponential type. 
This kind of problems will be the main focus of this paper.

\subsection{General scheme for construction of frames from kernels}
\label{construction-needlets}

The main application of the kernels $L_n(x,y)$ defined in (\ref{def:Ln-gen})
is to the construction of frames of rapidly decaying elements (needlets).
We now briefly describe the main elements of this construction.


\noindent
{\em $\bullet$ Semi-continuous Calder\'{o}n type decomposition.}
Suppose $\ha$ is an admissible cutoff function of type (c) in the sense of
Definition~\ref{cutoff-d1}. Then
$\sum_{\nu=0}^\infty |\ha(2^{-\nu}t)|^2= 1$, $t \in [1, \infty)$.
In the general setup of \S \ref{localizatin-principle}, define
$$
L_0(x,y) := P_0(x,y) \quad \mbox{and} \quad
L_j(x, y) := \sum_{\nu=0}^\infty
\wh a \Big(\frac{\nu}{2^{j-1}}\Big)P_\nu(x,y),
\quad  j=1, 2, \dots,
$$
and denote briefly
$
(L_j*f)(x):= \int_E L_j(x, y)f(y)d\mu(y).
$
The following 
decomposition follows readily from the conditions on $\ha$
\begin{equation}\label{Calderon}
f=\sum_{j=0}^\infty L_j*L_j*f
\quad \mbox{for $ f\in L^2(E, \mu)$.}
\end{equation}


\noindent
{\em $\bullet$ Discretization via cubature formulas.}
Suppose that there is a cubature formula
\begin{equation}\label{cubature-form}
\int_E f d\mu \sim \sum_{\xxi \in \cX_j} c_\xxi f(\xxi)
\end{equation}
with $\cX_j \subset E$ and $c_\xxi >0$,
which is exact for all functions $f$
of the form $f=gh$ with
$g, h \in \bigoplus_{\nu=0}^{2^{2j}} \CV_{\nu}$.
Cubature formula (\ref{cubature-form}) allows to rewrite (\ref{Calderon})
in the form
\begin{equation}\label{semi-discret-dec}
f(x)=\sum_{j=0}^\infty \sum_{\xi\in\cX_j}c_\xi^{1/2}L_j(\xi, x)
\int_Ef(y)c_\xi^{1/2}L_j(\xi, y)dy.
\end{equation}


\noindent
{\em $\bullet$ Definition of frame elements $($needlets$)$.}
Now the frame elements are defined by
\begin{equation}\label{def-needlets}
\psi_\xi(x):= c_\xi^{1/2}\cdot L_j(\xi, x)
\quad \mbox{for}\quad
\xi\in \cX_j, ~ j=0, 1, \dots.
\end{equation}
We write $\cX := \cup_{j = 0}^\infty \cX_j$, where equal points from different levels
$\cX_j$ are considered as distinct point of $\cX$, so that we can use $\cX$ as an index set
in the definition of the needlet system
$$
\Psi:=\{\psi_\xi\}_{\xi\in\cX}.
$$
From (\ref{semi-discret-dec}) and the definition of $\{\psi_\xi\}$ it readily follows that
$$
f= \sum_{\xi\in \cX} \langle f, \psi_\xi\rangle \psi_{\xi}
\quad\mbox{in}\;\; L^2(E, \mu)\quad\mbox{and}\quad
\|f\|_{L^2(E, \mu)}
=\Big(\sum_{\xi\in \cX} |\langle f, \psi_\xi\rangle|^2\Big)^{1/2},
$$
i.e. $\Psi$ is a tight frame for $L^2(E, \mu)$.

From definition (\ref{def-needlets}) it is clear that the frame elements $\psi_\xi$
inherit the rapid decay of the kernels $L_j$ if this is the case.
The superb localization of the building blocks $\{\psi_\xi\}$ is the reason for calling
them {\em needlets}.
The~rapid decay of needlets makes them a powerful tool for decomposition
of spaces of functions and distributions in various settings.
The above scheme has already been utilized for construction of needlets
and needlet decomposition of $L^p$, Sobolev, and the more general
Triebel-Lizorkin and Besov spaces in the frameworks of
spherical harmonics \cite{NPW1, NPW2},
Jacobi polynomials \cite{PX1, KPX1},
orthogonal polynomials on the ball \cite{PX2, KPX2},
and Hermite and Laguerre functions \cite{Epp2, PX3, KPPX}.

\subsection{Sub-exponentially localized wavelets, kernels, and needlets}
\label{sub-exponential-kern}

In \cite{DzH} Dziuba\'{n}ski and Hern\'{a}ndez constructed band-limited wavelets
of sub-exponential decay. More precisely, they showed that for any $\eps>0$
there exists a $C^\infty$ mother wavelet $\psi$ such that its Fourier transform
$\widehat{\psi}$ is compactly supported on $\RR$ and
\begin{equation}\label{local-wavelets}
|\psi(x)| \le c_\eps \exp\{-|x|^{1-\eps}\}, \quad x\in\RR,
\end{equation}
with $c_\eps$ a constant depending on $\eps$.
They also showed that in this estimate $\eps>0$ cannot be removed.
However, as will be shown in \S\ref{Examples}, the localization of $\psi$
can be improved to
$|\psi(x)| \le c \exp\{-\frac{c'|x|}{\ln (1+|x|)^{1+\eps}}\}$
and beyond (see \S\ref{Examples}).

Our aim in this paper is to construct kernels and needlets with similar
(sub-exponential) localization in the context of
Jacobi polynomials,
spherical harmonics,
orthogonal polynomials on the $d$-dimensional ball and simplex with weights,
and $d$-dimensional Hermite and Laguerre functions.

For instance, we shall show that for any $\eps>0$ there exists
an admissible cutoff function $\ha$ of type (a), (b) or (c) in the sense of Definition~\ref{cutoff-d1}
such that the kernels
$L_n(x, y)$ from (\ref{Chebyshev-Kernel-intr}) satisfy
(see Theorem~\ref{thm:Chebyshev-localization-1} below)
\begin{equation}\label{Sub-exp-loc-introd}
|L_n(x, y)|
\le c n \exp\left\{-\frac{c_\eps n \dd(x, y)}{[\ln(e+n \dd(x, y))]^{1+\eps}}\right\},
\quad x, y \in [-1, 1].
\end{equation}
Evidently, this estimate yields (\ref{Loc-Chebyshev-intr}).
The above estimate can be further improved by replacing the term
$[\ln(e+n \dd(x, y))]^{1+\eps}$
by any product of the form
\begin{align}\label{product2}
&\ln(e+ n\dd(x, y))\cdots
\underbrace{\ln\cdots\ln}_\ell \big(\underbrace{\exp\cdots\exp}_\ell 1+n\dd(x, y)\big)\\
&\hspace{1.4in}\times
\big[\underbrace{\ln\cdots\ln}_{\ell+1}\big(\underbrace{\exp\cdots\exp}_{\ell+1}1+n\dd(x, y)\big)\big]^{1+\eps},
\quad \ell>1.\notag
\end{align}
Estimate (\ref{Sub-exp-loc-introd}) leads to the following localization of
the $j$th level Chebyshev needlets
$
|\psi_\xi(x)|
\le c2^{j/2}\exp\left\{-\frac{c_\eps 2^j \dd(x, \xi)}{[\ln(e+2^j\dd(x, \xi))]^{1+\eps}}\right\},
$
which can be further improved as above.
We shall also show that the above estimates are sharp in the sense that
$\eps>0$ cannot be removed.

We would like to emphasize that according to the localization principle
kernels of the form e.g. (\ref{Chebyshev-Kernel-intr}) are rapidly decaying for
an arbitrary admissible cutoff function $\ha$, while the sub-exponential localization
of these kernels is only possible for exceptional cutoff functions~$\ha$.
One of the main steps in constructing kernels and frames of sub-exponential localization
is the construction of admissible cutoff functions $\ha$, most importantly ones of type (c),
with derivatives obeying
$\|\ha^{(k)}\|_\infty\le \cc(\cc/\eps)^kk^k(\ln k)^{k(1+\eps)}$,~$k\ge 3$,
see Theorem~\ref{thm:cutoff} below.
Although we employ various techniques for proving our results, the proof of
the sub-exponential localization of Jacobi kernels plays a~prominent role in this paper.


An outline of this paper is as follows:
%
In \S\ref{Examples} we clarify the impact of the behavior of the cutoff functions
around zero on the localization of the respective kernels in the simplest case
of Chebyshev polynomials.
We also illustrate the notion of sub-exponential localization of kernels
in the case of Chebyshev polynomials.
In \S\ref{not-valid-principle} we show that the localization principle described in
\S\ref{localizatin-principle} is no longer valid
in the case of tensor product Chebyshev or Legendre polynomials ($d=2$) and
products of Chebyshev and Legendre polynomials.
In \S\ref{cutoff} we construct admissible cutoff functions of ``small" derivatives.
Sections \ref{Jacobi-kernels}, \ref{sphere-kernels}, \ref{ball-kernels},
\ref{simplex-kernels}, \ref{Hermite-kernels}, and \ref{Laguerre-kernels} are devoted
to the construction of sub-exponentially localized kernels and needlets in the context of
Jacobi polynomials, spherical harmonics, orthogonal polynomials on the ball and simplex,
and Hermite and Laguerre functions.

Throughout this paper we shall use the following notation:
For $x\in\RR^d$ we shall use the norms
$\|x\|=\|x\|_\infty:=\max_i |x_i|$, $\|x\|_2:= (\sum_i |x_i|^2)^{1/2}$, and
$|x|=\|x\|_1:=\sum_i |x_i|$.
Positive constants will be denoted by $c$, $c_1$, $c'$, $\dots$ and they may vary
at every occurrence, $A\sim B$ will stand for $c_1A\le B\le c_2A$.
Also, $\lfloor x \rfloor$ will denote the largest integer not exceeding $x$.

\section{Localization principle and sub-exponential localization: Simple examples}
\label{Examples}
\setcounter{equation}{0}


We {\bf first} would like to clarify the impact of the behavior of
the cutoff function $\ha$ at $t=0$ on the localization properties of kernels
as in (\ref{def:Ln-gen}).
To this end, we shall use the simple example of normalized Chebyshev polynomials
of first kind on $[-1, 1]$ defined by
$$
\TT_n(x):= \sqrt{2/\pi}\cos n\arccos x
\quad \mbox{for $n\ge 1$ and}\quad
\TT_0:= \sqrt{1/\pi}.
$$
As is well-known $\{\TT_n\}_{n\ge 0}$ is an orthonormal basis for $L^2([-1, 1], (1-x^2)^{-1/2})$.
We are interested in the localization of the kernels $L_n(x, y)$,
defined in (\ref{Chebyshev-Kernel-intr})
%
%
with $\ha \in C^\infty[0, \infty)$ and, say, $\supp \, \ha \subset [0, 2]$.
Setting $x=:\cos \theta$ and $y=:\cos \phi$, we have
\begin{align*}
L_n(\cos \theta, \cos \phi)
= \frac{2}{\pi}\Big(\ha(0)/2 + \sum_{j=1}^\infty \ha\Big(\frac{j}{n}\Big) \cos\theta\cos\phi\Big)
= \frac{1}{\pi}(F_n(\theta-\phi)+F_n(\theta+\phi)),
\end{align*}
where
$
F_n(\theta):= \ha(0)/2 + \sum_{j=1}^\infty \ha(\frac{j}{n}) \cos\theta.
$
Let $\ha$ be the even extension of $\ha$, i.e. $\ha(t):= \ha(-t)$ for $t<0$.
Then
\begin{equation}\label{defn-Fn}
F_n(\theta)=\frac{1}{2}\sum_{j\in\ZZ}\ha\Big(\frac{j}{n}\Big)e^{ij\theta}
\end{equation}
and the question of localization of $L_n(x,y)$ reduces to the localization of
the trigonometric polynomial $F_n(\theta)$ around $\theta=0$.

It is easy to see (see the proof of Theorem~\ref{thm:Chebyshev-localization-1} below)
that if $\ha\in C^\infty$, then for any $\sigma>0$
there exists a constant $c_\sigma>0$ such that
\begin{equation}\label{localization-F}
|F_n(\theta)|\le c_\sigma\frac{n}{(1+n|\theta|)^\sigma},
\quad \theta\in [-\pi, \pi],
\end{equation}
which readily leads to
\begin{equation}\label{Chebyshev-localization}
|L_n(x,y)|\le c_\sigma\frac{n}{(1+n\dd(x, y))^\sigma}
\quad \forall \sigma>0,
\end{equation}
where $\dd(x, y)$ is the distance on $[-1, 1]$ defined by
$
\dd(x, y):= |\arccos x - \arccos y|.
$
The point is that the localization from (\ref{Chebyshev-localization})
is not possible if the even extension of $\ha$ to $\RR$ is not in $C^\infty$.
The precise argument is that (\ref{Chebyshev-localization}) holds for $y=1$
if and only if (\ref{localization-F}) is true.
In turn (\ref{localization-F}) is valid if and only if $\ha\in C^\infty$.
This last claim can be justified as follows:
Let $f$ be defined with Fourier transform
$\hat f(\xi):=\ha(\xi/n)e^{i\xi t}$.
Then
$f(y)= na\left(n(y+t)\right)$
and the Poisson summation formula
\begin{equation}\label{formula-poisson}
\sum_{j\in\ZZ} g(2\pi j) = (2\pi)^{-1}\sum_{j\in\ZZ}\hat g(j)
\quad \Big(\hat g(\xi):=\int_\RR g(t) e^{-i\xi t}dt\Big)
\end{equation}
gives
\begin{equation}\label{poisson-a}
F_n(t):= \frac{1}{2}\sum_{j\in\ZZ}
\ha\Big(\frac{j}{n}\Big)e^{ijt}
= \pi n\sum_{j\in\ZZ}a\left(n(t+2\pi j)\right).
\end{equation}
Hence, the inverse Fourier transform $a$ of $\ha$, has to be in the Schwartz class,
which is only possible if $\ha\in C^\infty$.
For this it suffices to have $\ha^{(k)}(0)=0$ for $k\ge 1$.
For other orthogonal systems (see below), however, it is not completely clear
what behavior of $\ha$ at $t=0$ would lead to rapid decay of the respective
kernels and needlets.

In concluding, it is evident that the behavior of the cutoff function $\ha$ at $t=0$
is important and not every compactly supported $C^\infty$ function $\ha$
will give rise to rapidly decaying kernels $L_n(x, y)$.
For our purposes, however, it suffices to restrict the selection of cutoff functions $\ha$ to
compactly supported $C^\infty$ functions which are constants in a neighborhood of $t=0$,
which automatically resolves the issue about the behavior of $\ha$ at $t=0$
(see Definition~\ref{cutoff-d1}).

\medskip


{\bf Secondly}, we want to illustrate with the next theorem the notion of
a sub-exponential localization of kernels of type (\ref{def:Ln-gen}) on the simple example
of Chebyshev polynomials of first kind.


\begin{theorem}\label{thm:Chebyshev-localization-1}
For any $0<\eps \le 1$ there exists an admissible cutoff function $\ha$
of type $(a)$ or $(b)$ or $(c)$
such that the kernels $L_n(x,y)$ from $(\ref{Chebyshev-Kernel-intr})$ satisfy
\begin{equation}\label{Chebishev-bound1}
|L_n(x,y)|
\le cn\exp\left\{-\frac{c'\eps n\dd(x, y)}{[\ln(e+n\dd(x, y))]^{1+\eps}}\right\},
\quad x, y\in [-1, 1],
\end{equation}
where
$c'>0$ is an absolute constant
and $c$ depends only on $\eps$.
\end{theorem}

\noindent
{\bf Proof.}
A key ingredient in this proof
is the existence of an admissible cutoff function $\ha$ of an arbitrary type
such that
\begin{equation}\label{derivative-a}
\|\ha^{(k)}\|_\infty\le \cc(\cc/\eps)^kk^k(\ln k)^{k(1+\eps)}
\quad \mbox{for $k\ge 3$},
\end{equation}
where $\cc >0$ is an absolute constant
(see Theorem~\ref{thm:cutoff} below).
Denote again by $\ha$ the even extension of $\ha$ on $\RR$.

Set $x=:\cos \theta$ and $y=:\cos \phi$.
Exactly as above
\begin{equation}\label{Repr-Ln}
L_n(\cos \theta, \cos \phi)
= \frac{1}{\pi}(F_n(\theta-\phi)+F_n(\theta+\phi)),
\end{equation}
where $F_n$ is defined in (\ref{defn-Fn}) and also (\ref{poisson-a}) holds.
Note that $F_n(\theta)$ is a $2\pi$-periodic trigonometric
polynomial of degree $< 2n$.
Evidently,
$$
t^\ell a(t)
= \frac{i^\ell}{2\pi} \int_\RR \ha^{(\ell)}(\xi)e^{i\xi t}d\xi,
$$
and hence
$|t|^\ell|a(t)|\le \|\hat a^{(\ell)}\|_\infty$
as $\supp \ha \subset [-2, 2]$.
This for $\ell=0$ and $\ell=k$ gives
$$
|a(t)| \le
2^k\frac{\|\ha^{(k)}\|_\infty+\|\ha\|_\infty}{(1+|t|)^k}
\le\frac{c(2\cc/\eps)^kk^k(\ln k)^{k(1+\eps)}}{(1+|t|)^k},
\quad k\ge 3,
$$
on account of (\ref{derivative-a}).
Combining this with (\ref{poisson-a}) gives
\begin{align}\label{est-Fn-first}
|F_n(\theta)|
&\le cn(2\cc/\eps)^kk^k(\ln k)^{k(1+\eps)}
\sum_{j\in\ZZ} \frac{1}{(1+n|\theta+2\pi j|)^k}\\
&\le cn\frac{(2\cc/\eps)^kk^k(\ln k)^{k(1+\eps)}}{(1+n|\theta|)^k}
= cn\left(\frac{(2\cc/\eps)k(\ln k)^{(1+\eps)}}{1+n|\theta|}\right)^k.\notag
\end{align}
We next use the above to show that
\begin{equation}\label{est-Fn}
|F_n(\theta)|
\le cn\exp\left\{-\frac{c'\eps n|\theta|}{[\ln(e+n|\theta|)]^{1+\eps}}\right\},
\quad |\theta|\le \pi.
\end{equation}
Indeed, if
$
1+n|\theta| \le 6e(2\cc/\eps)[\ln(e+n|\theta|)]^{1+\eps},
$
then (\ref{est-Fn}) follows from the obvious estimate
$|F_n(\theta)|\le 2n$.

Suppose
$
1+n|\theta| > 6e(2\cc/\eps)[\ln(e+n|\theta|)]^{1+\eps}.
$
Choosing
$k:=\Big\lfloor \frac{1+n|\theta|}{2e(2\cc/\eps)[\ln(e+n|\theta|)]^{1+\eps}}\Big\rfloor$
one easily shows that $k\ge 3$ and
$
\frac{(2\cc/\eps) k(\ln k)^{1+\eps}}{1+n|\theta|} \le e^{-1},
$
and (\ref{est-Fn}) follows by (\ref{est-Fn-first}).

Note that $\theta$, $\phi$ from (\ref{Repr-Ln}) obey $0\le \theta, \phi\le \pi$.
We use (\ref{est-Fn}) to estimate $|F_n(\theta-\phi)|$.

If $\theta+\phi\le \pi$, then we use (\ref{est-Fn}) and that
$\theta+\phi\ge |\theta-\phi|$ to estimate $|F_n(\theta+\phi)|$.

If $\theta+\phi > \pi$, then $0<2\pi-\theta-\phi\le \pi$
and since $F_n$ is $2\pi$-periodic and even, we have
$F_n(\theta+\phi)= F_n(2\pi-\theta-\phi)$.
We now use $2\pi-\theta-\phi \ge |\theta-\phi|$ and (\ref{est-Fn})
to estimate $|F_n(\theta+\phi)|=|F_n(2\pi-\theta-\phi)|$.
Putting together these estimates, we get
\begin{align*}
|L_n(\cos \theta, \cos \phi)|
&\le \frac{1}{\pi}(|F_n(\theta-\phi)|+|F_n(\theta+\phi))|)\\
&\le cn\exp\left\{-\frac{c'\eps n|\theta-\phi|}{[\ln(e+n|\theta-\phi|)]^{1+\eps}}\right\},
\end{align*}
which implies (\ref{Chebishev-bound1}).
$\qed$

As a byproduct of the above proof we get the following localization result for
trigonometric polynomials.
For a compactly supported cutoff function $\ha$ consider the trigonometric
polynomial
\begin{equation}\label{trig-polynomial}
L_n(\theta):=\sum_{j\in\ZZ}\ha\Big(\frac{j}{n}\Big)e^{ij\theta}.
\end{equation}


\begin{theorem}\label{thm:trig-localization-1}
For any $0<\eps \le 1$ there exists a $C^\infty$ cutoff function $\ha\not\equiv 0$
supported on $[-2,2]$
such that the polynomial $L_n(\theta)$ from $(\ref{trig-polynomial})$ satisfies
\begin{equation}\label{trig-bound1}
|L_n(\theta)|
\le cn\exp\left\{-\frac{c'\eps n|\theta|}{[\ln(e+n|\theta|)]^{1+\eps}}\right\},
\quad \theta\in[-\pi, \pi],
\end{equation}
where
$c'>0$ is an absolute constant
and $c$ depends only on $\eps$.
\end{theorem}

The construction of $\ha$ is as above and the proof of (\ref{trig-bound1})
is exactly the same as the proof of (\ref{est-Fn}) above.
We omit it.


Several remarks are in order.
Estimate $(\ref{trig-bound1})$ can be slightly improved,
namely, the term $[\ln (e+ n|\theta|)]^{1+\eps}$ in $(\ref{trig-bound1})$
can be replaced by
$\ln(e+ n|\theta|) [\ln\ln (e^e+n|\theta|)]^{1+\eps}$ or, in general, by any product of the form
\begin{equation}\label{product}
\ln(e+ n|\theta|)\cdots
\underbrace{\ln\cdots\ln}_\ell (\underbrace{\exp\cdots\exp}_\ell 1+n|\theta|)
\big(\underbrace{\ln\cdots\ln}_{\ell+1}
\big(\underbrace{\exp\cdots\exp}_{\ell+1}1+n|\theta|\big)\big)^{1+\eps}
\end{equation}
with $\ell \ge 1$.

Estimate $(\ref{Chebishev-bound1})$ can also be improved in a similar fashion
(see Theorem~\ref{thm:Jacobi-localization-1} below).
Note also that in these cases inequality (\ref{trig-bound1}) is trivially true for $\varepsilon=0$.

On the other hand, Theorem \ref{thm:trig-localization-1} is sharp in the sense that
(\ref{trig-bound1}) cannot be replaced by an estimate of this type
\begin{equation}\label{lower-bound}
|L_n(\theta)|\le c|L_n(0)|e^{-c'\ph(n|\theta|)},\quad \theta\in[-\pi,\pi],
\end{equation}
with $\ph$ of the form
$\ph(t)=t\ln(e+t)^{-1}$ or $\ph(t)=t\ln(e+t)^{-1}\ln\ln(e^e+t)^{-1}$ or etc.
This follows by Example 2.2 in \cite{Ivanov-Totik}, which implies that the inequality
\begin{equation}\label{lower-bound-2}
\int_1^\infty \frac{\ph(t)}{t^2}\,dt < \infty
\end{equation}
is a necessary condition for the validity of (\ref{lower-bound}) for an appropriate polynomial
$L_n$ of degree $n$.
A similar observation applies to Theorem \ref{thm:Chebyshev-localization-1} as well.

\medskip


Our {\bf third} point is that the sub-exponential localization of the wavelet
from (\ref{local-wavelets}), established in \cite{DzH}, can be
improved as follows.


\begin{theorem}\label{thm:wavelet-localization}
For any $0<\eps \le 1$ there exists a bandlimited orthogonal wavelet $\psi$
such that
\begin{equation}\label{wavelet-localiz}
|\psi(x)|
\le c\exp\left\{-\frac{c'\eps |x|}{[\ln(e+|x|)]^{1+\eps}}\right\},
\quad x \in \RR,
\end{equation}
where $c'>0$ is an absolute constant
and $c$ depends only on $\eps$.
Furthermore, the term $[\ln(e+|x|)]^{1+\eps}$ can be replaced by
a product similar to the one in $(\ref{product})$ above.
\end{theorem}

\noindent
{\bf Proof.} In essence, this proof is contained in the proof of
Theorem~\ref{thm:Chebyshev-localization-1}.
So, we only outline it here.
Choose $\ha\ge 0$ to be the admissible function of type (c) from the proof
Theorem~\ref{thm:cutoff} below
and such that its derivatives obey (\ref{derivative-a}).
Let us now denote again by $\ha$ the even extension of $\ha$ to all of $\RR$ and
set $\theta(\xi):=\ha(\frac{3}{4\pi}\xi)$.
From the properties of $\ha$ (see Theorem~\ref{thm:cutoff})
it readily follows that the function $\psi$ with Fourier transform
$$
\hat\psi(\xi):=\theta(\xi)e^{-i\xi/2}=\ha\Big(\frac{3}{4\pi}\xi\Big)e^{-i\xi/2}
$$
is a band limited orthogonal wavelet (see \cite{Meyer}).
It remains to show that $\psi$ obeys (\ref{wavelet-localiz}).
Indeed, we have
$$
\psi(x) = \frac{1}{2\pi}
\int_\RR \ha\Big(\frac{3}{4\pi}\xi\Big)e^{-i\xi/2} e^{ix\xi}d\xi,
$$
which implies
$$
(x-1/2)^\ell\psi(x)
= \frac{1}{2\pi}\Big(\frac{3i}{4\pi}\Big)^\ell
\int_\RR \ha^{(\ell)}\Big(\frac{3}{4\pi}\xi\Big)e^{i(x-1/2)\xi}d\xi
\quad \forall \ell.
$$
Similarly as in the proof of Theorem~\ref{thm:Chebyshev-localization-1},
this and (\ref{derivative-a}) yield
$$
|\psi(x)|
\le c\left(\frac{(2\cc/\eps)k(\ln k)^{(1+\eps)}}{1+|x|}\right)^k,
\quad k\ge 3,
$$
and one completes the proof by choosing $k$ appropriately.
\qed

\smallskip

Finally, we remark that estimate (\ref{wavelet-localiz}) is sharp in
the sense that it cannot be replaced by an estimate of the form
$
|\psi(x)|
\le c\exp\big\{-\frac{c'|x|}{\ln(e+|x|)}\big\}.
$
To see this, we shall use the argument we used to prove the sharpness
of estimate (\ref{trig-bound1}).
Precisely as in (\ref{poisson-a}) with $a(t)$ replaced by $\psi(x)$
(using the Poisson summation formula (\ref{formula-poisson})), we have
$$
\phi_n(x):=\sum_{j\in\ZZ}
\hat\psi\Big(\frac{j}{n}\Big)e^{ijx}
= 2\pi n\sum_{j\in\ZZ}\psi\left(n(x+2\pi j)\right).
$$
Note that since $\hat\psi$ is compactly supported, then $\phi_n$ from above
is a trigonometric polynomial of degree $O(n)$.
Assuming now that
$
|\psi(x)|
\le c\exp\big\{-\frac{c'|x|}{\ln(e+|x|)}\big\}
$
leads to
$
|\phi_n(x)|\le cn\exp\big\{-\frac{c'n|x|}{\ln(e+n|x|)}\big\}.
$
But this is impossible since the condition given in (\ref{lower-bound})-(\ref{lower-bound-2})
would be violated.

\section{Construction of $C^\infty$ cutoff functions}
\label{cutoff}
\setcounter{equation}{0}

The present section lays down some of the ground work that will be needed for
the construction of sup-exponentially localized kernels and needlets.
Admissible cutoff functions of ``small" derivatives will be needed.


\begin{theorem}\label{thm:cutoff}
Let $0<\eps\le 1$. Then there exists an admissible cutoff function
$\ha$ of any type: $(a)$, $(b)$ or $(c)$ $($see Definition~\ref{cutoff-d1}$)$
such that $0\le \ha\le 1$,
\begin{equation}\label{est-derivative}
\|\ha^{(k)}\|_\infty\le \cc(\cc/\eps)^kk^k(\ln k)^{k(1+\eps)}, \; k\ge 3,
\mbox{ and }\;
\|\ha^{(k)}\|_\infty\le \cc(\cc/\eps)^kk^k,
\; k=1,2.
\end{equation}
where $\cc>1$ is an absolute constant $($e.g. $\cc=88$$)$.
Moreover, $\ha^{(k)}(1)=0$ for $k\ge 1$.

Furthermore, there exists an admissible function $\ha$
of type $(a)$, $(b)$ or $(c)$
such that the above estimates still hold with the term
$(\ln k)^{k(1+\eps)}$ replaced by $(\ln k)^k(\ln\ln k)^{k(1+\eps)}$
or, in general, by a product of the form
$
(\ln k)^k \cdots
(\underbrace{\ln\cdots\ln}_\ell k)^k
(\underbrace{\ln\cdots\ln}_{\ell+1} k)^{k(1+\eps)}
$
with $\ell\ge 1$, for sufficiently large $k$ and
$\cc$ depending on $\ell$.
\end{theorem}

\noindent
{\bf Proof.}
We first construct a $C^\infty$ bump $h$ of ``small" derivatives.
To this end let
$\chi_\delta:= \frac{1}{2\delta}\ONE_{[-\delta, \delta]}$
and choose
$\delta_0=\delta_1:=1$ and
$\delta_j:=\frac{1}{j(\ln j)^{1+\eps}}$ if $j\ge 2$.
Note that
$$ 
\sum_{j=0}^\infty \delta_j\le 2+ \frac{1}{2}+\frac{1}{3} +
\int_3^\infty\frac{dt}{t(\ln t)^{1+\eps}}
< 3+\frac{1}{\eps(\ln 3)^\eps}
\le \frac{4}{\eps}.
$$ 
We define
$$
h_m:=\chi_{\delta_0}* \dots * \chi_{\delta_m}
\quad\mbox{and}\quad
h(t):=\lim_{m\to\infty} h_m(t).
$$
It is easy to see that (cf. \cite[Theorem 1.3.5]{H})
$h\in C^\infty$, $h\ge 0$, $\supp h\subset (-\frac{4}{\eps}, \frac{4}{\eps})$ and
$$
\|h^{(k)}\|_\infty \le \frac{1}{\prod_{j=0}^k \delta_j}
\le k^k(\ln k)^{k(1+\eps)}
\quad\mbox{for}\quad k\ge 3.
$$
Also, since $\int_\RR \chi_\delta=1$, we have $\int_\RR h=1$.

Our second step is to rescale $h$, namely, we define
$h_\eps(t):= \frac{8}{\eps}h(\frac{8 t}{\eps})$.

We next apply standard wavelet techniques. We first integrate $h_\eps$,
i.e. we define $g(t):=\frac{\pi}{2}\int_{-\infty}^th_\eps(s)ds$.
Evidently,
$g\in C^\infty$, $0\le g\le \pi/2$, $\supp g' \subset (-\frac{1}{2}, \frac{1}{2})$,
$g(t)+g(-t)=\frac{\pi}{2}$ for $t\in\RR$,
$\|g^{(k)}\|_\infty \le \frac{\pi}{2}(\frac{8}{\eps})^{k+1}$ for $0\le k\le 3$, and
\begin{equation}\label{est-gk}
\|g^{(k+1)}\|_\infty \le \frac{\pi}{2}\|h_\eps^{(k)}\|_\infty
\le \frac{\pi}{2}\Big(\frac{8}{\eps}\Big)^{k+1}k^k(\ln k)^{k(1+\eps)}
\quad\mbox{for}\quad k\ge 3.
\end{equation}
Observe that $\ha(t):=\frac{2}{\pi}g(\frac 3 2-t)$ is an admissible function of type $(a)$.

In~order to construct an admissible function $\ha$ of type $(c)$ and,
hence, of type $(b)$
we define $\phi(t):= \sin g(t)$, $t\in \RR$.
From above,
$\phi(t)^2+\phi(-t)^2=1$ for $t\in\RR$.
We now set
$$
\ha(t):=
\left\{
\begin{array}{lcl}
\phi(2t-\frac{3}{2}) &\mbox{if}& t\in [\frac{1}{2}, 1],\\
\phi(\frac{3}{2}-t) &\mbox{if}& t\in (1, 2],\\
0 &\mbox{if}& \RR\setminus [\frac{1}{2}, 2],
\end{array}
\right.
$$
and claim that $\ha$ has the properties of
an admissible function of type (c).
All of them are easy to verify but (\ref{est-derivative}),
which needs some care. Here we use some ideas from \cite{DzH}.
Fix $t_0\in (-\frac{1}{2}, \frac{1}{2})$ and let
$g_k(t):=\sum_{j=0}^k \frac{(t-t_0)^j}{j!}g^{(j)}(t_0)$ -
the $k$th degree Taylor polynomial of $g$ centered at~$t_0$.
Apparently, $\phi^{(k)}(t_0)= [\sin g_k]^{(k)}(t_0)$
and since $\sin g_k(z)$ is an entire function, by the Cauchy formula,
\begin{equation}\label{Cauchy}
\phi^{(k)}(t_0)=\frac{k!}{2\pi i}\int_\gamma \frac{\sin g_k(z)}{(z-t_0)^{k+1}}dz,
\end{equation}
where $\gamma:=\{z\in \bC:|z-t_0|=r\}$ with $r=\frac{1}{2e(8/\eps)(\ln k)^{1+\eps}}$.
By Stirling's formula $n!>(n/e)^n$
and using (\ref{est-gk}) we have for $z\in \gamma$ and $k\ge 3$
\begin{align*}
|g_k(z)|\le \sum_{j=0}^k \frac{1}{j!}\frac{\pi}{2}
\Big(\frac{8}{\eps}\Big)^{j}\frac{j^j(\ln k)^{j(1+\eps)}}{[2e(8/\eps)(\ln k)^{1+\eps}]^j}
\le \frac{\pi}{2}\sum_{j=0}^k\frac{1}{2^j}<\pi
\end{align*}
and, therefore,
$|\sin g_k(z)|\le e^\pi$ for $z\in \gamma$.
We use this in (\ref{Cauchy}) to obtain
$$
|\phi^{(k)}(t_0)| \le \frac{k!}{2\pi}e^\pi 2\pi [2e(8/\eps)(\ln k)^{1+\eps}]^k
\le e^\pi(44/\eps)^kk^k(\ln k)^{k(1+\eps)},
$$
which implies (\ref{est-derivative}).
Here we used that $k!\le k^k$.

For the construction of $\ha$ which satisfies (\ref{est-derivative}) with
$(\ln k)^{k(1+\eps)}$ replaced by a~product of the form
%
$
(\ln k)^k \cdots
(\underbrace{\ln\cdots\ln}_\ell k)^k
(\underbrace{\ln\cdots\ln}_{\ell+1} k)^{k(1+\eps)}
$
one proceeds as above with obvious modifications.
$\qed$


\begin{remark}
{\rm
Theorem \ref{thm:cutoff} cannot be improved in the sense that $\eps$ in the exponent
in $(\ref{est-derivative})$ cannot be removed $($replaced by $\eps=0$$)$.
Indeed, by the Denjoy-Carleman theorem $($see e.g. Theorem 1.3.8 in \cite{H}$)$
one has that if $\ha\in C^\infty$ on a given interval and
$$
\|\ha^{(k)}\|_\infty\le c^{k+1}M_k,
\quad k=0,1,\dots,
$$
then the divergence of the series
$\sum_{j=0}^\infty 1/L_j$ with $L_j=\inf_{k\ge j}M_k^{1/k}$ yields that
$\ha$ is quasi-analytic $($see \cite[Definition 1.3.7.]{H}$)$.
Having in mind that the cutoff functions are not quasi-analytic we conclude that
the above inequalities cannot be true with
$M_k=[k\ln(e+k)]^k$ or
$
M_k=\big[k\ln (e+ k)\cdots\underbrace{\ln\cdots\ln}_\ell
\big(\underbrace{\exp\cdots\exp}_\ell 1 + k)\big]^k,
$
in general.
}
\end{remark}

\section{Sub-exponentially localized kernels and frames induced by Jacobi polynomials}
\label{Jacobi-kernels}
\setcounter{equation}{0}

The Jacobi polynomials $\{P_n^{(\a,\b)}\}_{n\ge 0}$,
form an orthogonal basis for the weighted space
$L^2([-1, 1], w_{\a,\b})$ with weight $w_{\a,\b}(t):=(1-t)^\a(1+t)^\b$,
$\a, \b>-1$.
They are traditionally normalized by
$P_n^{(\a,\b)}(1)=\binom{n+\a}{n}$.
It is well known that~\cite[(4.3.3)]{Sz}
$$
\int_{-1}^1 P_n^{(\a,\b)}(t)
P_m^{(\a,\b)}(t)w_{\a,\b}(t)dt  = \delta_{n, m} h_n^{(\a,\b)},
$$
where
\begin{equation}\label{def-hn}
h_n^{(\a,\b)} =
\frac{2^{\a+\b+1}}{(2n+\a+\b+1)}
\frac{\Gamma(n+\a+1)\Gamma(n+\b+1)}{\Gamma(n+1)\Gamma(n+\a+\b+1)}.
\end{equation}

We are interested in kernels of the form
\begin{equation}\label{def.L}
L_n^{\a,\b}(x,y)=\sum_{j=0}^\infty \ha \Big(\frac{j}{n}\Big)
     \Big(h_j^{(\a,\b)}\Big)^{-1} P_j^{(\a,\b)}(x) P_j^{(\a,\b)}(y),
\end{equation}
where $\ha$ is an admissible cutoff function in the sense of Definition~\ref{cutoff-d1}.
In \cite{PX1} it is proved that this sort of kernels decay rapidly away from the main
diagonal in $[-1, 1]^d\times [-1, 1]^d$.
Our aim here is to show that for suitable admissible cutoff functions $\ha$
the localization of $L_n^{\a,\b}(x,y)$ can be improved to sub-exponential.

To state our result, we need the quantity
\begin{equation}\label{Jacobi-weight}
w_{\a,\b}(n; x)
:= (1-x + n^{-2})^{\a+1/2} (1+x + n^{-2})^{\b+1/2}.
\end{equation}
We shall also use again the distance
$\dd(x, y):=|\arccos x-\arccos y|$ on $[-1, 1]$.


\begin{theorem}\label{thm:Jacobi-localization-1}
Let $\alpha, \beta \ge -1/2$ and $0<\eps \le 1$.
Then for any admissible cutoff function $\ha$ obeying inequality $(\ref{est-derivative})$
in Theorem~\ref{thm:cutoff}
the kernels from $(\ref{def.L})$ satisfy the following inequality
for $x, y\in [-1, 1]$
\begin{equation}\label{Jacobi-bound1}
|L_n^{\a,\b} (x,y)|
\le \frac{cn}{\sqrt{w_{\a,\b}(n; x)}\sqrt{w_{\a,\b}(n; y)}}
\exp\Big\{-\frac{\CC n\dd(x, y)}{[\ln(e+n\dd(x, y))]^{1+\eps}}\Big\},
\end{equation}
where
$\CC=c'\eps$ with
$c'>0$ an absolute constant,
and $c$ depends only on $\a$, $\b$, and~$\eps$.

Furthermore, for an appropriate cutoff function $\ha$ the term
$[\ln(e+n\dd(x, y))]^{1+\eps}$ above can be replaced by any product of the form
$(\ref{product2})$.
\end{theorem}

For the proof of Theorem~\ref{thm:Jacobi-localization-1} and the respective
localization results on the sphere (Theorem~\ref{thm:sphere-localization-1})
and ball (Theorem~\ref{thm:ball-localization-1}), we first need to establish
the sub-exponential localization of $L_n^{\a,\b}(x,1)$.
Set
\begin{equation}\label{def.Ln1}
\Q_n^{\a,\b}(x) := L_n^{\a,\b}(x,1)
=\sum_{j=0}^\infty \ha\Big(\frac{j}{n}\Big)
  \Big(h_j^{(\a,\b)}\Big)^{-1} P_j^{(\a,\b)}(1)P_j^{(\a,\b)}(x).
\end{equation}
It is easy to see that
\begin{equation}\label{eq:L-n}
\Q_n^{\a,\b}(x) =  c^\star
\sum_{j=0}^\infty \ha\Big(\frac{j}{n}\Big)
\frac{(2j+\a+\b+1)\Gamma(j+\a+\b+1)}{\Gamma(j+\b+1)}P_j^{(\a,\b)}(x),
\end{equation}
where
$c^\star:=2^{-\a-\b-1}\Gamma(\a+1)^{-1}$.


\begin{theorem} \label{thm:Jacobi-localization-2}
Let $\a \ge \b \ge -1/2$ and $0<\eps \le 1$.
Assume that the cutoff function $\ha$ in $(\ref{def.Ln1})$ is given by Theorem~\ref{thm:cutoff}.
Then we have
\begin{equation}\label{est.Ln}
\Big|\frac{d^r}{dx^r} \Q_n^{\a,\b}(\cos \theta)\Big|
\le c n^{2\a+2r+2}\exp\left\{-\frac{\CC n\theta}{[\ln (e+ n\theta)]^{1+\eps}}\right\},
\quad 0 \le \theta \le \pi,
\end{equation}
where $\CC=c'\eps$ with $c'>0$ an absolute constant and
$c=c''8^r$ with $c''>0$ depending only on $\a$, $\b$, and $\eps$.

Moreover, for an appropriate cutoff function $\ha$ of type $(a)$, $(b)$, or $(c)$
from Theorem~\ref{thm:cutoff} the above estimate holds with the term
$[\ln (e+ n\theta)]^{1+\eps}$ replaced by any product of the form $(\ref{product})$.
\end{theorem}

\noindent
{\bf Proof.}
We first prove (\ref{est.Ln}) for $r=0$.
%
%
We obtain a trivial estimate by simply using that
$\Gamma(j + a) /\Gamma(j+1) \sim j^{a-1}$
and
$\|P_n^{(\a,\b)}\|_{L^\infty[-1, 1]} \le cn^\a$ with $c=c(\a, \b)$
\cite[(7.32.6)]{Sz}.
We get
\begin{equation}\label{est-Ln-triv}
|\Q_n^{\a,\b}(\cos\theta)|
     \le c \sum_{j=0}^{2n} j^{2 \a +1} \le c n^{2 \a+2},
\end{equation}
which implies (\ref{est.Ln}) ($r=0$) for $0\le \theta\le 1/n$.


A key role in proving a nontrivial estimate on $|\Q_n^{\a,\b}(\cos \theta)|$
will play the identity \cite[(4.5.3)]{Sz}:
\begin{align}\label{key-ident}
&\sum_{\nu=0}^n \frac{(2\nu+\a+k+\b+1)\Gamma(\nu+\a+k+\b+1)}{\Gamma(\nu+\b+1)}
P_\nu^{(\a+k,\b)}(x)\\
&\qquad\qquad= \frac{\Gamma(n+\a+k+1+\b+1)}{\Gamma(n+\b+1)}
P_n^{(\a+k+1,\b)}(x).\notag
\end{align}
Applying summation by parts to the sum in (\ref{eq:L-n})
(using (\ref{key-ident}) with $k=0$), we get
\begin{equation}\label{Ln-parts}
\Q_n^{\a,\b}(x) =  c^\star
\sum_{j=0}^\infty
\Big[\ha\Big(\frac{j}{n}\Big)-\ha\Big(\frac{j+1}{n}\Big)\Big]
\frac{\Gamma(j+\a+1+\b+1)}{\Gamma(j+\b+1)}P_j^{(\a+1,\b)}(x).
\end{equation}
We now define the sequence of functions $(A_k(t))_{k=0}^\infty$ by
$$
A_0(t):= (2t+\a+\b+1)\ha\Big(\frac{t}{n}\Big)
$$
and inductively
\begin{equation}\label{Ak+1}
A_{k+1}(t):= \frac{A_k(t)}{2t+\a+k+\b+1}-\frac{A_k(t+1)}{2t+\a+k+\b+3},
\quad k\ge 0.
\end{equation}
Notice that
\begin{equation}\label{A1}
A_1(t):= \ha\Big(\frac{t}{n}\Big)-\ha\Big(\frac{t+1}{n}\Big)
\end{equation}
and hence
$\supp A_k \subset [n/2-k, 2n] \subset [n/4, 2n]$ if $1\le k \le n/4$.

Applying summation by parts $k$ times starting from (\ref{eq:L-n})
(using every time (\ref{key-ident})), we arrive at the
identity:
\begin{equation}\label{Ln-parts-k}
\Q_n^{\a,\b}(x) = c^\star\sum_{j=0}^\infty
A_k(j)\frac{\Gamma(j+\a+k+\b+1)}{\Gamma(j+\b+1)}P_j^{(\a+k,\b)}(x).
\end{equation}


We next show that for all $1\le k \le n/4$ and $m \ge 0$
\begin{align}\label{est-Ak}
\|A_k^{(m)}\|_{\infty}
&\le \cc 12^{k-1}(e\cc/\eps)^{k+m}(k+m)^{1/2}(m+1)_{k-1}3^m m! n^{-2k-m+1}\\
&\hspace{2in}  \times [\ln (k+m+2)]^{(k+m)(1+\eps)}, \notag
\end{align}
where $\cc>1$ and $\eps$ are as in Theorem~\ref{thm:cutoff}
and $(a)_r:=a(a+1)\cdots(a+r-1)$.
In~particular,
\begin{align}\label{est-Ak0}
\|A_k\|_{\infty}
&\le \cc(12e\cc/\eps)^{k}k^{1/2}(k-1)![\ln (k+2)]^{k(1+\eps)}n^{-2k+1}\notag\\
&\le \cc(12e\cc/\eps)^{k}k^k[\ln (k+2)]^{k(1+\eps)}n^{-2k+1}
\quad\mbox{for}\quad 1\le k\le n/4.
\end{align}
To make the proof of (\ref{est-Ak}) more transparent we denote
$$
B_k^m:=\frac{n^{2k+m-1}}{3^mm!}\|A_k^{(m)}\|_{\infty}.
$$
Then (\ref{est-Ak}) is the same as
\begin{equation}\label{est-Bkm}
B_k^m
\le \cc 12^{k-1}(e\cc/\eps)^{k+m}(k+m)^{1/2}(m+1)_{k-1}
[\ln (k+m+2)]^{(k+m)(1+\eps)}.
\end{equation}
%
%
We shall proceed by induction on $k$.
By (\ref{A1}) and Theorem~\ref{thm:cutoff} it follows that
\begin{align*}
\|A_1^{(m)}\|_{\infty}
& \le n^{-m-1}\|\ha^{(m+1)}\|_\infty\\
&\le \cc(\cc/\eps)^{m+1}(m+1)^{m+1}[\ln (m+3)]^{(m+1)(1+\eps)}n^{-m-1}.
\end{align*}
But by Stirling's formula it follows that
$
(m+1)^{m+1} < (m+1)^{-1/2}e^{m+1}(m+1)!
$
and hence from above
$$
B_1^m
\le \cc(e\cc/\eps)^{m+1}(m+1)^{1/2}[\ln (m+3)]^{(m+1)(1+\eps)},
$$
which shows that (\ref{est-Bkm}) holds for $k=1$ and $m\ge 0$.

%
%
Suppose (\ref{est-Bkm}) holds for some $1\le k\le n/4$ and all $m\ge 0$.
We shall show that it holds with $k$ replaced by $k+1$ and for all $m\ge 0$.
Denote
$G_k(t):=\frac{A_k(t)}{2t+\a+k+\b+1}$.
Then by (\ref{Ak+1})
$
A_{k+1}^{(m)}(t)
=-\int_0^1G_k^{(m+1)}(t+s)ds
$
and evidently
$$
G_k^{(m+1)}(t)
=\sum_{\nu=0}^{m+1}\binom{m+1}{\nu}A_k^{(\nu)}(t)
\frac{(-2)^{m+1-\nu}(m+1-\nu)!}{(2t+\a+k+\b+1)^{m+2-\nu}}.
$$
Therefore,
\begin{align*}
\|A_{k+1}^{(m)}\|_\infty
\le \|G_k^{(m+1)}\|_\infty
&\le \sum_{\nu=0}^{m+1}\binom{m+1}{\nu}\|A_k^{(\nu)}\|_\infty
\frac{2^{m+1-\nu}(m+1-\nu)!}{(n-k)^{m+2-\nu}}\\
&\le 2\sum_{\nu=0}^{m+1}
\frac{(m+1)!3^{m+1-\nu}}{\nu!n^{m+2-\nu}}
\|A_k^{(\nu)}\|_\infty,
\end{align*}
where we used that since $k\le n/4$ we have $n-k> 3n/4$.
It is readily seen that the above is the same as
$$
B_{k+1}^m  \le 6(m+1)\sum_{\nu=0}^{m+1}B_k^\nu.
$$
We now use the inductive assumption (see   (\ref{est-Bkm}))
to obtain the following upper bound for $B_{k+1}^m$:
\begin{align*}
&6(m+1)\sum_{\nu=0}^{m+1}
 \cc 12^{k-1}(e\cc/\eps)^{k+\nu}(k+\nu)^{1/2}(\nu+1)_{k-1}
  [\ln (k+\nu+2)]^{(k+\nu)(1+\eps)}\\
&\le 6 \cc 12^{k-1}(k+m+1)^{1/2}(m+1)_k
  [\ln (k+m+3)]^{(k+m+1)(1+\eps)}
  \sum_{\nu=0}^{m+1}(e\cc/\eps)^{k+\nu}
\end{align*}
For the last sum above we have
$$
 \sum_{\nu=0}^{m+1}(e\cc/\eps)^{k+\nu}
 < \frac{(e\cc/\eps)^{k+m+2}}{e\cc/\eps-1}
 < 2(e\cc/\eps)^{k+m+1}
$$
using that $e\cc/\eps>2$.
Putting the above estimates together we get
$$
B_{k+1}^m
\le  \cc 12^k(e\cc/\eps)^{k+m+1}(k+m+1)^{1/2}(m+1)_k
  [\ln (k+m+3)]^{(k+m+1)(1+\eps)},
$$
which shows that (\ref{est-Bkm}) holds when $k$ is replaced by $k+1$.
Therefore, (\ref{est-Bkm}) and hence (\ref{est-Ak}) hold for all $1\le k\le n/4$ and $m\ge 0$.


We shall need the following estimate for Jacobi polynomials
\cite[Theorem 1]{EMN}: For $\a, \b \ge -1/2$ and $n\ge 1$,
\begin{equation}\label{est.Pn}
\sup_{x\in[-1, 1]}(1-x)^{\a+1/2}(1+x)^{\b+1/2}|P_n^{(\a,\b)}(x)|^2
\le \frac{2e}{\pi}\big(2+\sqrt{\a^2+\b^2}\big)h_n^{(\a,\b)},
\end{equation}
where $h_n^{(\a,\b)}$ is from (\ref{def-hn}).


We next prove (\ref{est.Ln}) ($r=0$) for $1/n \le \theta\le \pi/2$.
By (\ref{def-hn}) it readily follows that
$h_n^{(\a+k, \b)}\le c2^k/n$.
Using this, (\ref{est.Pn}), and the obvious inequality
$1-\cos \theta = 2\sin^2\theta/2\ge \theta^2$ for $0\le \theta\le \pi/2$,
we infer
\begin{equation}\label{est.Pn-k}
|P_n^{(\a+k,\b)}(\cos \theta)|
\le \frac{ck^{1/2}2^{k/2}}{n^{1/2}\theta^{k+\a+1/2}}
\le \frac{c2^k}{n^{1/2}\theta^{k+\a+1/2}},
\quad 0<\theta\le \pi/2.
\end{equation}
We now use this and (\ref{est-Ak0}) in (\ref{Ln-parts-k}) to obtain
for $1/n \le \theta \le \pi/2$ and $1\le k\le n/4$
\begin{align}\label{est-Ln-main}
|\Q_n^{\a,\b}(\cos \theta)|
&\le c(12e\cc/\eps)^{k}k^k[\ln (k+2)]^{k(1+\eps)}n^{-2k+1}
\sum_{n/2-k<j<2n} \frac{2^k(j+k)^{\a+k}}{j^{1/2}\theta^{k+\a+1/2}}\notag\\
&\le c(54e\cc/\eps)^{k}k^{k}[\ln (k+2)]^{k(1+\eps)}
\frac{n^{2\a+2}}{(n\theta)^{k+\a+1/2}}\\
&\le cn^{2\a+2}
\Big(\frac{c_\eps k[\ln (k+2)]^{1+\eps}}{n\theta}\Big)^k,
\quad c_\eps:= 54e\cc/\eps.\notag
\end{align}
Here we used that
$\Gamma(j+\a+k+\b+1)/\Gamma(j+\b+1)\le c(j+k)^{\a+k}$.

If $n\theta< 2ec_\eps[\ln(e+n\theta)]^{1+\eps}$, then (\ref{est.Ln}) (in the case $r=0$)
follows by (\ref{est-Ln-triv}) with $\CC=(2ec_\eps)^{-1}$.
Assume that $n\theta\ge 2ec_\eps[\ln(e+n\theta)]^{1+\eps}$
and choose
$k:=\Big\lfloor \frac{n\theta}{ec_\eps[\ln(e+n\theta)]^{1+\eps}}\Big\rfloor$.
Evidently, $k\le n/4$.
We claim that
\begin{equation}\label{optimization}
\frac{c_\eps k[\ln (k+2)]^{1+\eps}}{n\theta} \le e^{-1},
\end{equation}
Indeed, by the definition of $k$ we have
$k\le\frac{n\theta}{ec_\eps[\ln(e+n\theta)]^{1+\eps}}$
and then (\ref{optimization}) follows by the obvious inequality
$k+2 \le e+n\theta$.
Combining (\ref{optimization}) with (\ref{est-Ln-main}) leads to (\ref{est.Ln})
for $r=0$ with $\CC=(ec_\eps)^{-1}$.
Therefore, (\ref{est.Ln}) ($r=0$) holds  for $1/n \le \theta\le \pi/2$
with $\CC=(2ec_\eps)^{-1}=c'\eps$.


Let now $\pi/2 < \theta\le \pi-1/n$.
Similarly as in the proof of estimate (\ref{est.Pn-k}) we use that
$
1+\cos \theta=2\sin^2\frac{\pi-\theta}{2} \ge (\pi-\theta)^2
$
for $\pi/2\le \theta\le \pi$ to obtain
$$
|P_n^{(\a+k, \b)}(\cos\theta)|
\le \frac{c2^k}{n^{1/2}(\pi-\theta)^{\b+1/2}}
\le c2^kn^\beta,
\quad \pi/2\le \theta\le \pi-1/n.
$$
Combining this with (\ref{Ln-parts-k}) and (\ref{est-Ak0}) we get
\begin{align*}
|\Q_n^{\a,\b}(\cos \theta)|
&\le c(12e\cc/\eps)^{k}k^k[\ln (k+2)]^{k(1+\eps)}n^{-2k+1}
\sum_{n/2-k<j<2n}2^kn^\beta(j+k)^{\a+k}\\
&\le c(54e\cc/\eps)^{k}k^{k}[\ln (k+2)]^{k(1+\eps)}n^{-k+\a+\b+2}\\
&\le cn^{\a+\b+2}
\Big(\frac{c_\eps k[\ln (k+2)]^{1+\eps}}{n}\Big)^k,
\quad c_\eps:= 54e\cc/\eps.
\end{align*}
Exactly as above this leads to the estimate
\begin{equation}\label{est.Qn-3}
\Big|\Q_n^{\a,\b}(\cos \theta)\Big|
\le c n^{\a+\b+2}\exp\left\{-\frac{\CC n}{[\ln (e+ n)]^{1+\eps}}\right\},
\quad \pi/2 \le \theta \le \pi-1/n,
\end{equation}
and, since $\a\ge \b$,  estimate (\ref{est.Ln}) ($r=0$) holds in this case as well.


Finally, let $\pi-1/n\le \theta\le \pi$.
In this case  (\ref{est.Ln}) follows easily from  (\ref{est.Qn-3}).
Indeed, as is well known for any polynomial $P\in\Pi_n$ one has
$$
|P(x)|\le \|P\|_{L^\infty[-1, 1]}T_n(x)
\le  \|P\|_{L^\infty[-1, 1]}\big(x+\sqrt{x^2-1}\big)^n,
\quad x\ge 1,
$$
where $T_n$ is the $n$th degree Chebyshev polynomial of first kind,
and hence
$$
|P(x)|\le \|P\|_{L^\infty[-1, 1]} \big(1+2\sqrt{x-1}\big)^n
\quad \mbox{for}\quad 1\le x\le 5/4.
$$
By changing variables we obtain
\begin{equation}\label{est.P<Tn}
|P(t)|\le \|P\|_{L^\infty[a, b]}\Big(1+2^{3/2}\sqrt{(t-b)/(b-a)}\Big)^n,
\quad b\le t \le b+(b-a)/8.
\end{equation}
If $t:=-\cos \theta$, then $t\in [0, \cos \frac 1n]$ for $\theta\in [\pi-1/n, \pi]$.
On the other hand, $\cos \frac 1n = 1-2\sin^2 \frac{1}{n}\ge 1-2/n^2$.
Hence
\begin{equation}\label{Qn<Qn}
\|\Q_n^{\a,\b}(-\cdot)\|_{L^\infty[0, 1-2/n^2]}
\le \|\Q_n^{\a,\b}(\cos \cdot)\|_{L^\infty[\pi/2, \pi-1/n]}.
\end{equation}
Since $\Q_n^{\a,\b}\in\Pi_{2n}$, then we can apply (\ref{est.P<Tn}) on the interval
$[a, b]=[0, 1-2/n^2]$ to obtain using (\ref{Qn<Qn})
\begin{align*}
&\|\Q_n^{\a,\b}(\cos \cdot)\|_{L^\infty[\pi-1/n, \pi]}
\le \|\Q_n^{\a,\b}(-\cdot)\|_{L^\infty[1-2/n^2, 1]}\\
&\quad\le \|\Q_n^{\a,\b}(-\cdot)\|_{L^\infty[0, 1-2/n^2]}\
\Big(1+2^{3/2}\sqrt{(2/n^2)/(1-2/n^2) }\Big)^{2n}\\
&\quad \le (1+2^{5/2}/n)^{2n} \|\Q_n^{\a,\b}(\cos \cdot)\|_{L^\infty[\pi/2, \pi-1/n]}
\le c\|\Q_n^{\a,\b}(\cos \cdot)\|_{L^\infty[\pi/2, \pi-1/n]},
\end{align*}
provided $2/n^2 \le (1/8)(1-2/n^2)$ which is the same as $n\ge 5$
(the case $n<5$ is trivial).
Combining the above with (\ref{est.Qn-3}) implies that  the estimate in (\ref{est.Qn-3})
holds for $\pi-1/n\le \theta\le \pi$ as well, but perhaps with a different constant $c$.
This completes the proof of estimate (\ref{est.Ln}) when $r=0$.


For $r\ge 1$, estimate (\ref{est.Ln}) is an easy consequence
of Markov's inequality: If $Q \in \Pi_m$, then
$
\|Q'\|_{L^\infty[a,b]} \le 2m^2(b-a)^{-1}\|Q\|_{L^\infty[a,b]}.
$
Indeed, substituting $x=\cos \theta$ shows that (\ref{est.Ln}) with $r=0$
is the same as
\begin{equation}\label{last-est-Ln}
|\Q_n^{\a,\b}(x)|
\le \frac{c n^{2\a+2}}{E(\arccos x)}, 
\quad\mbox{where}\;\;
E(\theta):=\exp\left\{
\frac{\CC n\theta}{[\ln (e+ n\theta)]^{1+\eps}}
\right\}.
\end{equation}
It is easy to see that the function $E(\theta)$ is increasing on $[0, \pi]$
and hence $1/E(\arccos x)$ is also increasing on $[-1, 1]$.
Note that $\Q_n^{\a,\b}\in\Pi_{2n}$ and by Markov's inequality
we have, for $x \in [0,1]$,
\begin{align*}
\Big|\frac{d}{dx}\Q_n^{\a,\b}(x)\Big|
\le 8n^2(1+x)^{-1} \|\Q_n^{\a,\b}\|_{L^\infty[-1,x]}
\le 8n^2 \frac{c n^{2\a+2}}{E(\arccos x)},
\end{align*}
where we used (\ref{last-est-Ln}) and the monotonicity of $1/E(\arccos x)$.
Iterating we get
\begin{align*}
\Big|\frac{d^r}{dx^r}\Q_n^{\a,\b}(x)\Big|
\le c8^rn^{2r}\frac{n^{2\a+2}}{E(\arccos x)},
\quad r\ge 1,
\end{align*}
which yields (\ref{est.Ln}) for $0\le \theta\le \pi/2$.

To estimate $|(d/dx)^r \Q_n^{\a,\b}(x)|$ for $x \in [-1,0)$
one applies Markov's inequality on $[-1, 0]$ which leads readily to the same estimate
with a different constant $\CC$.

We finally observe that one proves estimate (\ref{est.Ln}) with
$[\ln (e+ n\theta)]^{1+\eps}$ replaced by a product of the form (\ref{product})
exactly as above using the respective estimate for $|\ha^{(k)}(t)|$ from Theorem~\ref{thm:cutoff}.
The proof of Theorem~\ref{thm:Jacobi-localization-2} is complete.
$\qed$

\medskip


\noindent
{\bf Proof of Theorem~\ref{thm:Jacobi-localization-1}.}
We shall use the notation and scheme of the proof of Theorem~2 in \cite{PX1}.
Consider first the case when $\a>\b>-1/2$.
By the product formula of Jacobi polynomials \cite{Koon} ($\a > \b > -1/2$), one has
\begin{equation} \label{product-Jac}
\frac{P_n^{(\a,\b)}(x)  P_n^{(\a,\b)}(y)}{ P_n^{(\a,\b)}(1)}
 = c_{\a,\b} \int_{0}^\pi \int_0^1  P_n^{(\a,\b)} (t(x,y,r,\psi)) dm(r,\phi),
\end{equation}
where
$$
 t(x,y,r,\psi) = \tfrac{1}{2}(1+x)(1+y)+ \tfrac{1}{2} (1-x)(1-y) r^2 +
      r\sqrt{1-x^2}\sqrt{1-y^2}  \cos \psi -1,
$$
the integral is against
$$
  d m(r,\phi) = (1-r^2)^{\a-\b-1} r^{2\b+1} (\sin \psi)^{2\b} dr d\psi,
$$
and the constant $c_{\a,\b}$ is selected so that
$
c_{\a,\b}\int_0^\pi\int_0^1 1 \,d m(r,\phi)=1.
$
Therefore, as in \cite[(2.15)]{PX1} we have
\begin{equation} \label{eq:Ln}
L_n^{\a,\b}(x, y) = c_{\a,\b} \int_0^\pi \int_0^1
\Q_n^{\a,\b}(t(x,y,r,\psi)) d m(r,\psi).
\end{equation}
Evidently, if $t = \cos \theta$ with $0 \le \theta \le \pi$, then
$\theta  \sim \sin \theta /2 \sim \sqrt{1-t}$.
Consequently, estimate (\ref{est.Ln}) implies
\begin{equation} \label{est-Qn-ab}
\left|\Q_n^{\a,\b}(t) \right|
\le cn^{2 \a +2}
\exp\left\{-\frac{c'\eps n\sqrt{1-t}}{[\ln(e+n\sqrt{1-t})]^{1+\eps}}
\right\}, \quad -1 \le t \le 1.
\end{equation}
Precisely as for the proof of (2.18) in \cite{PX1},
this yields
\begin{equation} \label{est-Qn1}
\left|L_n^{\a,\b}(\cos \theta, \cos \phi)\right|
\le cn^{2 \a +2}
\exp\left\{-\frac{c'\eps n|\theta-\phi|}{[\ln (e+n|\theta-\phi|)]^{1+\eps}}\right\},
\;\; 0\le \theta, \phi\le \pi.
\end{equation}
But
$
\exp\big\{\frac{-\frac{1}{2}c'\eps t}{[\ln(e+t)]^{1+\eps}}\big\}
\le c(1+t)^{-3\a-3\b-2},
$
$t>0$, for a sufficiently large constant $c>0$,
and hence (\ref{est-Qn1}) implies
\begin{equation} \label{est-Ln-ab}
\left|L_n^{\a,\b}(\cos \theta, \cos \phi)\right|
\le \frac{cn^{2 \a +2}}{(1+n|\theta-\phi|)^{3\a+3\b+2}}
\exp\left\{-\frac{\frac{1}{2}c'\eps n|\theta-\phi|}{[\ln (e+n|\theta-\phi|)]^{1+\eps}}\right\}.
\end{equation}
Just as in the proof of estimate (2.14) in \cite{PX1}
we use this estimate to prove (\ref{Jacobi-bound1}) with $\CC=c'/2$.
We skip the further details.

In the case when $\a=\b=-1/2$ (the case of Chebyshev polynomials of first kind),
estimate (\ref{Jacobi-bound1}) was already
proved in Theorem~\ref{thm:Chebyshev-localization-1}.

In the case when $\a=\b>-1/2$, as in \cite{PX1} and above, one uses
the product formula for Gegenbauer polynomials and (\ref{est-Ln-ab})
to prove estimate (\ref{Jacobi-bound1}). We omit the details.

Finally, let $\a>\b=-1/2$.
Evidently, for a continuous function $f$ on $[-1, 1]$
$$
  \lim_{\mu \to -1^+} \frac{\int_{0}^\pi  f(\cos \theta) (\sin \theta)^\mu d\theta}
   {\int_0^\pi (\sin \theta)^\mu d\theta} =  \frac{f(1) + f(-1)}{2}.
$$
On the other hand,
$
\lim_{\beta\to -1/2} P_n^{(\a,\b)}(x)=P_n^{(\a,-1/2)}(x)
$
uniformly on every compact interval, and in particular on $[-1, 1]$,
since as is well know the coefficients of $P_n^{(\a,\b)}$ converge
to the respective coefficients of $P_n^{(\a,-1/2)}$ as $\beta\to -1/2$.
Therefore, passing to the limit in (\ref{product-Jac}) as $\beta\to -1/2$
the product formula of Jacobi polynomials takes the form
\begin{align}\label{product-1/2}
&\frac{P_n^{(\a,-1/2)}(x)  P_n^{(\a,-1/2)}(y)}{ P_n^{(\a,-1/2)}(1)}\notag\\
& \qquad \qquad \qquad
=\frac{1}{2}c_{\a} \int_0^1  \left[ P_n^{(\a,-1/2)} (t(x,y,r,0)) + P_n^{(\a,-1/2)} (t(x,y,r, \pi))\right]
    dm(r),
\end{align}
where $d m(r) =  (1-r^2)^{\a-1/2}dr$ and $c_{\a}\int_0^1 \,d m(r)=1$.
Using this product formula, one can carry over the entire proof of Theorem 2 in [20]
in the case when $\alpha > \beta = -1/2$.
In a similar fashion, one uses (\ref{est-Ln-ab})-(\ref{product-1/2})
to prove estimate (\ref{Jacobi-bound1}) whenever $\alpha > \beta = -1/2$.
We skip the details.

Due to symmetry, in all other cases for $\a$, $\b$ estimate (\ref{Jacobi-bound1})
follows from the cases considered above.
\qed

\smallskip

\noindent
{\bf Sub-exponentially localized needlets induced by Jacobi polynomials.}
Given $j\ge 0$ let $\cX_j$ be the set of all zeros of the Jacobi polynomial
$P^{(\a, \b)}_m$ of degree $m:=2^j$ and
let $\{c_\xi\}_{\xi\in\cX_j}$ be the coefficients of the Gaussian quadrature formula
on $[-1, 1]$, which is exact for all polynomials of degree $2m-1$ (cf. (\ref{cubature-form})).
Suppose $\ha$ is an admissible cutoff function of type (c) obeying $(\ref{est-derivative})$.
According to the general scheme from \S\ref{construction-needlets} the $j$th level needlets
are defined by
$$
\psi_\xi(x):=c_\xi^{1/2}L_n^{\a,\b}(\xi, x),
\quad \xi\in\cX_j,
$$
where $L_n^{\a,\b}$ with $n:=2^{j-1}$ is the kernel defined in (\ref{def.L}).
The needlet system $\Psi:=\{\psi_\xi\}$ is a tight frame for $L^2([-1, 1], w_{\a, \b})$
and it is easy to see that the $j$th level needlet inherit from the kernel
$L_n^{\a,\b}$ the localization:
\begin{equation}\label{local-Needlet-Jacobi}
|\psi_\xi(x)|
\le \frac{c 2^{j/2}}{\sqrt{w_{\a,\b}(2^j; \xi)}}
\exp\Big\{-\frac{c'\eps 2^j\dd(x, \xi)}{[\ln(e+2^j\dd(x, \xi))]^{1+\eps}}\Big\},
\quad \xi\in\cX_j.
\end{equation}
Here we used the notation from Theorem~\ref{thm:Jacobi-localization-1}.
This is an improvement compared with the localization of
the Jacobi needlets constructed in \cite{PX1} and \cite{KPX1}.

\section{Sub-exponentially localized kernels and frames on the sphere}
\label{sphere-kernels}
\setcounter{equation}{0}

Denote by $\cH_n$ the space of all spherical harmonics of degree $n$
on the $d$-dimensional unit sphere $\SS^d$ in $\RR^{d+1}$.
As is well known (see e.g. \cite{SW}) the kernel of the orthogonal projector onto $\cH_n$
is given by
\begin{equation}\label{def-Pnu}
\PP_n(\xi\cdot\eta)= \frac{n+\lambda}{\lambda\omega_d}
C_n^\lambda (\xi\cdot \eta),
\end{equation}
where $\lambda:=\frac{d-1}{2}$,
$\omega_d:=\int_{\SS^d} 1 d\sigma$
is the hypersurface area of $\SS^d$,
and $\xi\cdot\eta$ stands for the inner product of $\xi, \eta\in \SS^d$.
Here $C_n^\lambda$ is the Gegenbauer polynomial of degree
$n$ normalized with
$C_n^\lambda(1)= \binom{\nu + 2\lambda-1}{\nu}$
\cite[p.~174]{Edelyi_v2}.

Our aim is to construct sub-exponentially localized kernels of the form
\begin{equation}\label{def-LN}
L_n(\xi\cdot\eta) = \sum_{j=0}^\infty \ha\Big(\frac{j}{n}\Big)\PP_j(\xi\cdot\eta),
\end{equation}
where $\ha$ is a $C^\infty$ cutoff function.

Denote by $\dd(\xi, \eta):=\arccos (\xi\cdot\eta)$ the geodetic distance between
$\xi, \eta\in \SS^d$ on $\SS^d$.


\begin{theorem}\label{thm:sphere-localization-1}
Let $0<\eps \le 1$.
Then for any admissible cutoff function $\ha$ obeying inequality $(\ref{est-derivative})$
in Theorem~\ref{thm:cutoff} the kernels from $(\ref{def-LN})$ satisfy
\begin{equation} \label{est-local-sphere}
\left|L_n(\xi\cdot\eta)\right|
\le cn^d
\exp\left\{-\frac{\CC n \dd(\xi, \eta)}{[\ln(e+n \dd(\xi, \eta))]^{1+\eps}}\right\},
\quad \xi, \eta\in\SS^d,
\end{equation}
where $\CC=c'\eps$ with
$c'>0$ an absolute constant
and $c>0$ depends only on $d$ and~$\eps$.
Moreover, for an appropriate cutoff function $\ha$
the above estimate can be improved similarly as in Theorem~\ref{thm:Jacobi-localization-1}.
\end{theorem}


\noindent
{\bf Proof.}
By (\ref{eq:L-n}) with $\a=\b=\lambda-1/2$ ($\lambda:=(d-1)/2$) we get
$$
\Q_n^{\lambda-1/2,\lambda-1/2}(t) =  c^\star
\sum_{j=0}^\infty \ha\Big(\frac{j}{n}\Big)
\frac{(j+\lambda)\Gamma(j+2\lambda)}{\Gamma(j+\lambda+1/2)}P_j^{(\lambda-1/2,\lambda-1/2)}(t),
$$
where
$c^\star:=2^{-2\lambda+1}\Gamma(\lambda+1/2)^{-1}$.
On the other hand, using the relation between Gegenbauer and Jacobi polynomials
\cite[(4.7.1)]{Sz} we have
\begin{equation}\label{relation-C-P}
C_n^\lambda(t)
=\frac{\Gamma(\lambda+1/2)}{\Gamma(2\lambda)}
 \frac{\Gamma(n+2\lambda)}{\Gamma(n+\lambda+1/2)}
P_n^{(\lambda-1/2, \lambda-1/2)}(t).
\end{equation}
Therefore,
$
L_n(t)=c(d)\Q_n^{\lambda-1/2,\lambda-1/2}(t)
$
and hence estimate (\ref{est-local-sphere}) follows immediately from
Theorem~\ref{thm:Jacobi-localization-2}.
$\qed$

\smallskip

\noindent
{\bf Sub-exponentially localized needlets on the sphere.}
As in \cite{NPW2} for the construction of the $j$th level needlets on $\SS^d$
we use a cubature formula with nodes in $\cX_j\subset \SS^d$ consisting of $O(2^{jd})$
almost uniformly distributed points on the sphere and positive coefficients
$\{c_\xi\}_{\xi\in\cX_j}$ of size $c_\xi\sim 2^{-jd}$,
which is exact for all spherical harmonics of degree $\le 2^{j+1}$
(for more details, see \cite{NPW2}).
Then choosing an admissible cutoff function $\ha$ of type (c)
satisfying (\ref{est-derivative}), we define
$$
\psi_\xi(x):=c_\xi^{1/2}L_n(\xi\cdot x),
\quad \xi\in\cX_j,
$$
where $L_n$ with $n:=2^{j-1}$ is the kernel defined in (\ref{def-LN}).
Setting $\cX:=\cup_{j\ge 0}\cX_j$, we define the needlet system by
$\Psi:=\{\psi_\xi\}_{\xi\in\cX}$.
This is a tight frame for $L^2(\SS^d)$ and Theorem~\ref{thm:sphere-localization-1}
implies the sub-exponential localization of the spherical needlets:
\begin{equation} \label{needled-local-sphere}
\left|\psi_\xi(x)\right|
\le c2^{jd/2}
\exp\left\{-\frac{c'\eps 2^j\dd(x, \xi)}{[\ln(e+2^j\dd(x, \xi))]^{1+\eps}}\right\},
\quad x\in\SS^d,\quad \xi\in\cX_j.
\end{equation}
Here, $\dd(x, \xi)$ is the geodesic distance between $x, \xi\in\SS^d$.
This is a natural improvement of the localization of the needlets
constructed in \cite{NPW2}.

\section{Sub-exponentially localized kernels and frames on the ball}
\label{ball-kernels}
\setcounter{equation}{0}

Here, we consider orthogonal polynomials on the unit ball
$B^d:=\{x\in\RR^d: \|x\|_2<1\}$ in $\RR^d$ $(d>1)$
with weight
$$
w_\mu(x):= (1-\|x\|_2^2)^{\mu-1/2}, \quad \mu \ge 0.
$$
We shall use the notation and results from \cite{PX2}.
Denote by $\CV_n^d$ the space of all polynomials of total degree $n$ which
are orthogonal to lower degree polynomials in $L^2(B^d, w_\mu)$.
As is shown in \cite{X99}, if $\mu>0$ the orthogonal projector
$\proj_n: L^2(B^d, w_\mu) \mapsto\CV_n^d$
can be written in the form
\begin{equation}\label{Proj-Vn}
(\proj_n f)(x) = \int_{B^d} f(y)P_n(w_\mu; x,y) w_\mu(y) dy,
\end{equation}
where
\begin{align*}\label{eq:Kn}
P_n(w_\mu; x,y) = c_{\mu, d}
\frac{\lambda+n}{\lambda}
\int_{-1}^1 C_n^\lambda \left(\langle x,y\rangle +
u\sqrt{1-\|x\|_2^2} \sqrt{1-\|y\|_2^2}\right) (1-u^2)^{\mu-1}du.
\end{align*}
Here $\langle x, y \rangle$ is the Euclidean inner product in $\RR^d$,
$C_n^\lambda$ is the $n$th degree Gegenbauer polynomial,
and
$$
\lambda = \mu + \frac{d-1}{2}.
$$
In the case $\mu = 0$ the kernel $P_n(w_\mu; x,y)$ is simpler, see \cite{PX2}.

As before, we are interested in the construction of sup-exponentially localized
kernels of the form
\begin{equation}\label{def:LnBd}
L_n^\mu(x,y)=\sum_{j=0}^\infty \ha\Big(\frac{j}{n}\Big)P_j(w_\mu;x,y),
\quad x,y \in B^d,
\end{equation}
where $\ha$ is a $C^\infty$ cutoff function.
As in \cite{PX2}, we shall need the following distance on $B^d$:
\begin{equation}\label{dist-Bd}
\dd(x,y):= \arccos \left
  \{ \langle x,y\rangle + \sqrt{1-\|x\|_2^2}\sqrt{1-\|y\|_2^2} \right \}
\end{equation}
and the quantity
\begin{equation*}\label{def.W}
W_\mu(n; x) := \left(\sqrt{1-\|x\|_2^2} + n^{-1}\right)^{2\mu}.
\end{equation*}


\begin{theorem}\label{thm:ball-localization-1}
Let $0<\eps \le 1$.
Then for any admissible cutoff function $\ha$ obeying inequality $(\ref{est-derivative})$
in Theorem~\ref{thm:cutoff} the kernels from $(\ref{def:LnBd})$ satisfy
\begin{equation}\label{ball-bound1}
|L_n^\mu (x,y)|
\le \frac{cn^d}{\sqrt{W_\mu(n; x)}\sqrt{W_\mu(n; y)}}
\exp\Big\{-\frac{\CC n\dd(x, y)}{[\ln(e+n\dd(x, y))]^{1+\eps}}\Big\},
\end{equation}
where
$\CC=c'\eps$ with $c'>0$ an absolute constant,
and $c$ depends only on $\mu$, $d$, and~$\eps$.

Moreover, for an appropriate cutoff function $\ha$
the term $[\ln(e+n\dd(x, y))]^{1+\eps}$ above can be replaced by any product of the form
$(\ref{product2})$.
\end{theorem}


\noindent
{\bf Proof.}
Assume $\mu>0$ (the case $\mu=0$ is easer).
From (\ref{Proj-Vn}) it follows that
\begin{equation} \label{rep.LnBd1}
L_n^\mu(x,y) = c(\mu, d)
\int_{-1}^1 \Q_n^\lambda (\langle x,y\rangle +
u\sqrt{1-\|x\|_2^2} \sqrt{1-\|y\|_2^2}) (1-u^2)^{\mu-1}du,
\end{equation}
where
$
\Q_n^\lambda(t) := \sum_{j=0}^\infty \ha\left(\frac{j}{n}\right)
          \frac{j+\lambda}{\lambda} C_j^\lambda (t).
$
Combining identities (\ref{relation-C-P}) and (\ref{eq:L-n}) gives
$
\Q_n^{\lambda-1/2, \lambda-1/2}(t)=c\Q_n^\lambda(t).
$
Therefore, by (\ref{est-Qn-ab}) we get
\begin{equation} \label{est-Qn-lambda}
\left|\Q_n^\lambda(t) \right|
\le cn^{2\lambda+1}
\exp\left\{-\frac{c'\eps n\sqrt{1-t}}{[\ln(e+n\sqrt{1-t})]^{1+\eps}}
\right\}, \quad -1 \le t \le 1,
\end{equation}
and since
$\exp\big\{\frac{-(c'/2)\eps u}{\ln (e+u)^{1+\eps}}\big\}\le c(1+u)^{-3\mu-1}$
for sufficiently large constant $c>0$
\begin{equation} \label{est-Qn-lambda1}
\left|\Q_n^\lambda(t) \right|
\le \frac{cn^{2\lambda+1}}{(1+n\sqrt{1-t})^{3\mu+1}}
\exp\left\{-\frac{\frac{1}{2}c'\eps n\sqrt{1-t}}{[\ln(e+n\sqrt{1-t})]^{1+\eps}}
\right\}.
\end{equation}
Now, just as in the proof of Theorem~4.2 in \cite{PX2}, we use estimates
(\ref{est-Qn-lambda})-(\ref{est-Qn-lambda1}) to obtain (\ref{ball-bound1}).
$\qed$

\smallskip

\noindent
{\bf Sub-exponentially localized needlets on the ball.}
For any $j\ge 0$ we shall utilize the cubature formula on $B^d$ from \cite{PX2}
with nodes in $\cX_j\subset B^d$
consisting of $O(2^{jd})$ almost uniformly distributed points on $B^d$ with respect
to the distance $\dd(\cdot, \cdot)$ defined in (\ref{dist-Bd})
and with positive coefficients $\{c_\xi\}_{\xi\in\cX_j}$, which is exact
for algebraic polynomials of degree $2^{j+1}$.
The second step is to select an admissible cutoff function $\ha$ of type (c)
as in Theorem~\ref{thm:ball-localization-1}.
Then we define
$$
\psi_\xi(x):=c_\xi^{1/2}L_n^\mu(\xi, x),
\quad \xi\in\cX_j,
$$
where $L_n^\mu$ with $n:=2^{j-1}$ is the kernel from (\ref{def:LnBd}).
Set $\cX:=\cup_{j\ge 0}\cX_j$.
Then we define the needlet system on the ball $B^d$ by
$\Psi:=\{\psi_\xi\}_{\xi\in\cX}$.
The sub-exponential localization of the needlets on the ball follows by
Theorem~\ref{thm:ball-localization-1}:
\begin{equation}\label{local-Needlet-ball}
|\psi_\xi(x)|
\le \frac{c 2^{j/2}}{\sqrt{W_\mu(2^j; \xi)}}
\exp\Big\{-\frac{c'\eps 2^j\dd(x, \xi)}{[\ln(e+2^j\dd(x, \xi))]^{1+\eps}}\Big\},
\quad x\in B^d, \quad \xi\in\cX_j.
\end{equation}
This is an improvement in comparison with the localization of the needlets
in~\cite{PX3}.

\section{Sub-exponentially localized kernels and frames on the simplex}
\label{simplex-kernels}
\setcounter{equation}{0}

Here, we consider orthogonal polynomials on the simplex
$$
T^d = \{x \in \RR^d: x_1 \ge 0, \ldots, x_d\ge 0, 1-|x|_1 \ge 0\},  \qquad
 |x|_1 := x_1 + \ldots + x_d,
$$
in $\RR^d$ ($d \ge 1$) with weight function
$$
w_{\k}(x) = x_1^{\k_1-\frac12} \cdots x_d^{\k_d-\frac12}(1-|x|_1)^{\k_{d+1}-\frac12}, \qquad
    \k_i \ge 0.
$$
Denote by $\CV_n^d$ the space of all polynomials of total degree $n$ that
are orthogonal to lower degree polynomials in $L^2(T^d,w_\k)$. As is shown in
\cite{X98} the orthogonal projector $\proj_n: L^2(T^d, w_\kappa) \mapsto \CV_n^d$
can be written in the form
\begin{equation}\label{Proj-VnT}
  (\proj_n f)(x) = \int_{T^d} f(y)P_n(w_\kappa; x,y) w_\kappa(y) dy,
\end{equation}
where, if all $\k_i > 0$,
\begin{align}\label{eq:KnT}
P_n(w_\kappa; x,y) = c(\kappa, d) \frac{2n+\lambda_\k}{\lambda_\k}
  \int_{[-1,1]^{d+1}} C_{2n}^{\lambda_\k} \left(z(x,y,t) \right) \prod_{i=1}^{d+1}
     (1-t_i^2)^{\k_i-1}d t.
\end{align}
Here
\begin{align}\label{eq:z(x,y,t)}
    z(x,y,t):= \sqrt{x_1y_1} t_1 + \cdots + \sqrt{x_dy_d} t_d +
       \sqrt{x_{d+1}y_{d+1}} t_{d+1}
\end{align}
with $x_{d+1} := 1-\|x\|_1$  and $y_{d+1} := 1-\|y\|_1$,
$C_n^\lambda$ is the $n$th degree Gegenbauer polynomial,
and
$$
       \lambda_\k := |\k| + \frac{d-1}{2}, \qquad |\k| : =\k_1+\cdots + \k_{d+1}.
$$
If some $\k_i=0$, then the identity in \eqref{eq:KnT} holds
with the integral in $t_i$ replaced according to the limit relation (see \cite{X98})
$$
\lim_{\kk \to 0} \int_{-1}^1 f(t) (1-t^2)^{\kk-1} dt \Big/\int_{-1}^1 (1-t^2)^{\kk-1} dt
= \frac{f(1)+ f(-1)}2.
$$

We are interested in kernels of the form
\begin{equation}\label{def:LnTd}
L_n^\k(x,y)=\sum_{j=0}^\infty \ha\Big(\frac{j}{n}\Big)P_j(w_\k;x,y),
\quad x,y \in T^d,
\end{equation}
where $\ha$ is a $C^\infty$ cutoff function.
Our aim is to show that for any admissible cutoff function $\ha$ the kernels $L_n^\k(x,y)$
decay rapidly away from the main diagonal in $T^d\times T^d$
and for any cutoff function $\ha$ from Theorem~\ref{thm:cutoff}
the decay is sub-exponential.

In analogy with the kernels on the ball, we shall need the distance $\ddT(\cdot, \cdot)$
on $T^d$ defined by
\begin{equation*}\label{def.distT}
\ddT(x,y):= \arccos \left
  \{ \sqrt{x_1 y_1} + \cdots + \sqrt{x_d y_d}+ \sqrt{x_{d+1}y_{d+1}} \right \}
\end{equation*}
and the quantity
\begin{equation*}\label{def.WT}
W_\k(n; x) :=  \prod_{i=1}^{d+1}  (x_i+n^{-2})^{\k_i}, \qquad x_{d+1} := 1-\|x\|_1.
\end{equation*}

Given $x, y\in T^d$, set
$a_j := \sqrt{x_j} =: \cos \theta_j$ and $b_j := \sqrt{y_j} =: \cos \phi_j$,
where $0\le \theta_j, \phi_j\le \pi/2$.
Applying the Cauchy-Schwartz inequality we get, for $1\le j\le d$,
$$
  \sum_{i=1}^d a_i b_i + \sqrt{1- a_1^2 -\cdots -a_d^2}\sqrt{1- b_1^2 -\cdots -b_d^2}
    \le a_j b_j + \sqrt{1-a_j^2} \sqrt{1-b_j^2}.
$$
Hence,
$$
\ddT(x,y) \ge  \arccos (\cos \theta_j \cos \phi_j + \sin \theta_j \sin \phi_j) =
      \arccos (\cos(\theta_j - \phi_j)),
$$
which yields $\ddT(x,y) \ge |\theta_j - \phi_j|$
and as a consequence
\begin{equation} \label{eq:dT}
   |\sqrt{x_j} - \sqrt{y_j}| \le \ddT(x,y), \quad 1\le j\le d+1.
\end{equation}

Our localization results take the form:


\begin{theorem}\label{thm:simplex-localization-0}
Let $\ha$ be an admissible cutoff function according to Definition~\ref{cutoff-d1}.
Then for any $\sigma>0$ there exists a constant $c_\sigma$ depending on
$\sigma$, $|\k|$, and $d$ such that the kernels from $(\ref{def:LnTd})$ satisfy

\begin{equation}\label{simplex-bound-0}
|L_n^\k (x,y)|
\le \frac{c_\sigma n^d}{\sqrt{W_\k(n; x)}\sqrt{W_\k(n; y)}}
\big(1 + n\ddT(x, y)\big)^{-\sigma} ,
\quad x,y\in T^d.
\end{equation}
\end{theorem}

The next theorem show that for suitable cutoff functions $\ha$ the localization of
the kernels $L_n^\k (x,y)$ can be improved to sub-exponential.


\begin{theorem}\label{thm:simplex-localization-1}
Suppose $0<\eps \le 1$ and let $\ha $ in $(\ref{def:LnTd})$ be an admissible cutoff function
obeying inequality $(\ref{est-derivative})$ in Theorem~\ref{thm:cutoff}.
Then 
\begin{equation}\label{simplex-bound1}
|L_n^\k (x,y)|
\le \frac{c n^d}{\sqrt{W_\k(n; x)}\sqrt{W_\k(n; y)}}
\exp\Big\{-\frac{\CC n \ddT(x, y)}{[\ln(e+n\ddT(x, y))]^{1+\eps}}\Big\},
\end{equation}
where
$\CC=c'\eps$ with $c'>0$ an absolute constant, and $c$ depends only
on $\k$, $d$, and $\eps$.

As elsewhere in this article, for an appropriate cutoff function $\ha$
estimate $(\ref{simplex-bound1})$ can be improved by replacing
the term $[\ln(e+n\ddT(x, y))]^{1+\eps}$ by any product of logarithms as in $(\ref{product2})$.
\end{theorem}

We shall only present the proof of Theorem~\ref{thm:simplex-localization-1}
since the proof of Theorem~\ref{thm:simplex-localization-0} follows along the same lines
but is simpler.

\medskip

\noindent
{\bf Proof of Theorem~\ref{thm:simplex-localization-1}.}
For a given $0<\eps\le 1$, let $\ha$ be the cutoff function from Theorem~\ref{thm:cutoff}.
Consider the case when $\kappa_i>0$ for $1\le i\le d+1$
(the case when some $\kappa_i$'s are zeros is treated in the same way with
appropriate modifications).
We begin with the relation
$$
  C_{2n}^\lambda(x) = \frac{\Gamma(n+ \lambda) \Gamma(\frac12)}
       {\Gamma(\lambda) \Gamma(n+\frac12)}
          P_n^{(\lambda -\frac12, -\frac12)} (2 x^2 -1),
$$
which follows readily combining identities (4.1.5) and (4.7.1) in \cite{Sz}.
This allows us to express $L_n^\k(x,y)$ in terms of the univariate kernel $\Q_n^{\a,\b}$
from (\ref{def.Ln1})-(\ref{eq:L-n}), namely,
\begin{equation*}
    L_n^\k(x,y) = c(\k,d) \int_{[-1,1]^{d+1}}
          \Q_n^{\lambda_\k -\frac12, -\frac12} \left(2 z(x,y,t)^2 -1\right)
              \prod_{i=1}^{d+1} (1-t_i^2)^{\k_i-1} dt.
\end{equation*}
Let $\theta(x,y,t) := \arccos (2z(x,y,t)^2 -1)$ with $z(x,y,t)$ given by \eqref{eq:z(x,y,t)}.
We use the fact that
$$
  1- z(x,y,t)^2  = \frac12 (1 - \cos \theta(x,y,t)) =  \sin^2 \frac{\theta(x,y,t)}{2}
       \sim \theta(x,y,t)^2
$$
and apply estimate \eqref{est-Qn-ab} for $\Q_n^{\a,\b}$ with $\alpha=\k-1/2$, $\beta=-1/2$
to obtain
\begin{equation*}
    |L_n^\k(x,y)| \le c n^{2 \la_\k+1} \int_{[-1,1]^{d+1}}
      \exp\left\{\frac{- c' \varepsilon n \sqrt{1-z(x,y,t)^2}}
        {[\ln (e+ n\sqrt{1-z(x,y,t)^2})]^{1+\eps}}\right\}
              \prod_{i=1}^{d+1} (1-t_i^2)^{\k_i-1} dt.
\end{equation*}
%
%
Evidently,
$$
1-z(x,y,t)^2 \ge 1-|z(x,y,t)|
\ge  1 - \sqrt{x_1y_1} |t_1| - \ldots - \sqrt{x_{d+1}y_{d+1} }|t_{d+1}|
$$
and using the symmetry of the integrand with respect to $t\in [-1, 1]^{d+1}$, we get
\begin{equation*}
 |L_n^\k(x,y)| \le c n^{2 \la_\k+1} \int_{[0,1]^{d+1}}
 \exp\left\{\frac{- c' \varepsilon n \sqrt{1-\zz(x,y,t)}}{\big[\ln \big(e+ n\sqrt{1-\zz(x,y,t)}\big)\big]^{1+\eps}}\right\}
              \prod_{i=1}^{d+1} (1-t_i^2)^{\k_i-1} dt.
\end{equation*}
Furthermore, we have the lower bound estimate
\begin{align*}
1- \zz(x,y,t)  & \ge 1- \sqrt{x_1y_1} - \ldots - \sqrt{x_{d+1}y_{d+1} }= 1- \cos \ddT(x,y)  \\
   &  = 2 \sin^2 \frac{\ddT(x,y)}{2} \ge \frac{2}{\pi^2} \ddT(x,y)^2, \notag
\end{align*}
which enables us to deduce the estimate
\begin{align} \label{eq:LnT.est}
&|L_n^\k(x,y)| \le  c n^{2 \la_\k+1}
    \exp\left\{-\frac{\frac{1}{2}c'\varepsilon n \ddT(x,y)}{[\ln (e+ n \ddT(x,y))]^{1+\eps}}\right\} \\
& \times \frac{1}{(1+n\ddT(x, y))^{|\k|}}
\int_{[0,1]^{d+1}} \frac{1}{ [ 1+ n \sqrt{1-\zz(x,y,t)} ]^\g}
\prod_{i=1}^{d+1} (1-t_i^2)^{\k_i-1} dt,    \notag
\end{align}
where $\gamma = 2|\k| + d +1$.
Here we used that
$\exp\big\{\frac{-(c'/2)\eps u}{\ln (e+u)^{1+\eps}}\big\}\le c(1+u)^{-|\k|-\gamma}$,
$u\ge 0$,
for sufficiently large constant $c>0$.
From the definition of $\zz(x,y,t)$, we have
$$
  1- \zz(x,y,t) = 1- \cos \ddT(x,y) + \sum_{i=1}^{d+1} \sqrt{x_iy_i} (1- t_i)
     \ge \sum_{i=1}^{d+1} \sqrt{x_iy_i} (1- t_i).
$$
Denote by $J$ the integrals in (\ref{eq:LnT.est}). From above, we have
\begin{align*}
J \le\int_{[0,1]^{d+1}} \frac{\prod_{i=1}^{d+1} (1-t_i^2)^{\k_i-1} dt}
{\left( 1+ n\big[\sum_{i=1}^{d+1} \sqrt{x_iy_i} (1- t_i)\big]^{1/2}\right)^\g}
=: I_{d+1}(\g)
\end{align*}
and our next goal is to estimate $I_{d+1}(\g)$.
To this end we first establish the following inequality for $A>0$, $B \ge 0$,
$\gamma \ge 2\k+1$, $\k>0$:
\begin{equation} \label{est-integral-T}
   \int_{0}^1 \frac{(1-t^2)^{\k -1} dt } {(1+ n \sqrt{B + A (1-t) })^\g} \le
       \frac{ cn^ {-2\k} } {A^\k \big(1+ n \sqrt{B}\big)^{\g - 2 \k -1}}.
\end{equation}
Indeed, substituting $s = n^2 A(1-t)$, we see that
\begin{align*}
&\int_{0}^1 \frac{(1-t^2)^{\k -1} dt } {\big(1+ n \sqrt{B + A(1-t) }\big)^\g}
\le \frac{2^{\k-1}}{(An^2)^\k} \int_{0}^{An^2} \frac{s^{\k -1} ds } {(1+  \sqrt{n^2 B + s })^\g} \\
& \le \frac{2^{\k-1}}{(An^2)^\k (1+ n\sqrt{B})^{\g-2 \k -1} }
      \int_{0}^{\infty} \frac{s^{\k -1} ds } {(1+  \sqrt{s })^{2\k +1}}
\le \frac{c n^{-2\k}} {A^\k \big(1+ n\sqrt{B}\big)^{\g-2 \k -1} }.
\end{align*}
We now set $B := 1 + n \sum_{i=1}^d \sqrt{x_iy_i} t_i$ and $A := \sqrt{x_{d+1}y_{d+1}}$,
and apply inequality \eqref{est-integral-T} to the integral in $I_{d+1}(\g)$
with respect to $t_{d+1}$. We get
\begin{align*}
I_{d+1}(\g)
& \le   \frac{ cn^ {-2\k_{d+1}} } {(\sqrt{x_{d+1}y_{d+1}})^{\k_{d+1}}}
  \int_{[0,1]^d} \frac{\prod_{i=1}^{d+1} (1-t_i^2)^{\k_i-1} dt}
       {\Big( 1+ n \big[\sum_{i=1}^d \sqrt{x_iy_i}(1-t_i)\big]^{1/2} \Big)^{\g-2\k_{d+1}-1} } \\
& =   \frac{c n^ {-2\k_{d+1}} } {(\sqrt{x_{d+1}y_{d+1}})^{\k_{d+1}}}
      I_d(\g - 2\k_{d+1} -1).
\end{align*}
Iterating this we obtain
$$
  I_{d+1}(\g) \le  \frac{cn^ {-2 |\k| } }
      {\prod_{i=1}^{d+1} (\sqrt{x_iy_i })^{\k_i}}.
$$
One the other hand, we trivially have $I_{d+1}(\g) \le 1$
and hence
$$
  I_{d+1}(\g) \le  \frac{cn^ {-2 |\k| } }
      {\prod_{i=1}^{d+1} (\sqrt{x_iy_i }+n^{-2})^{\k_i}}.
$$
Therefore,
\begin{equation}\label{est-J}
\frac{J}{(1+n\ddT(x, y))^{|\k|}}
\le  \frac{cn^ {-2 |\k| } }
      {(1+n\ddT(x, y))^{|\k|}\prod_{i=1}^{d+1} (\sqrt{x_iy_i }+n^{-2})^{\k_i}}.
\end{equation}
The simple inequality 
$$
  (a+n^{-1})(b+n^{-1}) \le 3 (ab+n^{-2})(1+n|a-b|), \quad a, b \ge 0, \quad n \ge 1,
$$
and \eqref{eq:dT} imply
\begin{align*}
\sqrt{x_i+n^{-2}} \sqrt{y_i+n^{-2}} &\le (\sqrt{x_i}+n^{-1})(\sqrt{y_i}+n^{-1})\\
\le 3(\sqrt{x_iy_i}+n^{-2})(1+n|\sqrt{x_i}-\sqrt{y_i}|)
&\le 3(\sqrt{x_iy_i}+n^{-2})(1+n\ddT(x, y)).
\end{align*}
This coupled with (\ref{est-J}) gives
$$
\frac{J}{(1+n\ddT(x, y))^{|\k|}}
\le \frac{cn^ {-2 |\k| } }{\sqrt{W_\k(n; x)}\sqrt{W_\k(n; y)}}.
$$
and inserting the above in (\ref{eq:LnT.est}) leads to (\ref{simplex-bound1}).
\qed

\smallskip

\noindent
{\bf Sub-exponentially localized needlets on the simplex.}
Needlet systems on the simplex $T^d$ have not been developed yet.
To this end one should follow the well established scheme from e.g. \cite {NPW2, PX1, PX2, PX3}.
The main ingredient of this development are
Theorems~\ref{thm:simplex-localization-0} and \ref{thm:simplex-localization-1}
from above, which provide the needed localization results.
Since some other elements of this theory are not completely developed yet we shall not
speculate here and leave the subject open.

\section{Sub-exponentially localized kernels and frames induced by Hermite functions}
\label{Hermite-kernels}
\setcounter{equation}{0}

The Hermite polynomials are defined by
$$
H_n(t) = (-1)^n e^{t^2} \frac{d^n}{dt^n} \Big(e^{-t^2}\Big),
\qquad n = 0, 1, \ldots .
$$
These polynomials are orthogonal on $\RR$ with weight $e^{-t^2}$.
We will denote by $h_n$ the $L^2$-normalized Hermite functions, i.e.
$$
h_n(t) := \left(2^n n! \sqrt{\pi}\right)^{-1/2} H_n(t) e^{-t^2/2}.
$$
Then the $d$-dimensional Hermite functions $\CH_\alpha(x)$ are defined by
$$
\CH_\alpha(x) := h_{\alpha_1}(x_1) \cdots h_{\alpha_d}(x_d),
\quad x=(x_1, \dots, x_d),
$$
where $\alpha=(\alpha_1, \dots, \alpha_d)\in \NN_0^d$.
Denoting by $\proj_n$ the orthogonal projector 
onto $\CV_n :={\rm span}\,\{\CH_\alpha: |\alpha|=n\}$, we have
$$
(\proj_n f)(x) = \int_{\RR^d} f(y) \CH_n(x,y)dy \quad \mbox{with}\quad
\CH_n(x,y):= \sum_{|\alpha| = n} \CH_\alpha(x)\CH_\alpha(y).
$$

It was shown in \cite{PX3} that for admissible cutoff functions $\ha$
the kernels
\begin{equation}\label{Hermire-kernels}
L_n(x,y):= \sum_{j=0}^\infty \ha\Big(\frac{j}{n}\Big) \CH_j(x,y)  
\end{equation}
decay rapidly away from the main diagonal in $\RR^d\times\RR^d$.
Our goal here is to show that for a suitable admissible cutoff functions $\ha$
the localization of $L_n(x, y)$ can be improved to sup-exponential.


\begin{theorem}\label{thm:Hermite-localization-1}
Let $0<\eps \le 1$ and assume that $\ha$ in $(\ref{Hermire-kernels})$
is an admissible cutoff function obeying inequality $(\ref{est-derivative})$
in Theorem~\ref{thm:cutoff}.
Then
\begin{equation}\label{Hermite-bound1}
|L_n (x,y)|
\le  c n^{d/2} \exp\Big\{-\frac{\CC n^{1/2} \|x-y\|} {[\ln(e+n^{1/2}\|x-y\| )]^{1+\eps}}\Big\},
\quad x,y\in  \RR^d,
\end{equation}
where
$\CC=c'\eps$ with $c'>0$ an absolute constant, and $c$ depends only
on $d$ and $\eps$. Recall that $\|x\|:=\max_j |x_j|$.

As before, for a suitable cutoff function $\ha$ the term
$[\ln(e+n^{1/2}\|x-y\| )]^{1+\eps}$ above
can be replaced by any product of multiple logarithmic terms as in $(\ref{product2})$.

Furthermore, for some constant $c''>0$
\begin{equation}\label{Hermite-bound2}
|L_n(x,y)|
\le ce^{-c''\max\{\|x\|^2, \|y\|^2\}}
\quad \mbox{if}\;\; \max\{\|x\|, \|y\|\}\ge (8n+2)^{1/2}.
\end{equation}

\end{theorem}

\noindent
{\bf Proof.}
As is well known (see \cite[Lemma 1.5.1]{Th}) for some constants $\gamma, c>0$
$$
|h_n(t)| \le e^{-\gamma t^2}
\quad\mbox{for}\quad t\ge (4n+2)^{1/2}
\quad\mbox{and}\quad
\|h_n\|_\infty \le cn^{-1/12}\le c,
$$
which readily yields (\ref{Hermite-bound2}).
In turn (\ref{Hermite-bound2}) immediately implies (\ref{Hermite-bound1})
in the case when $\max\{\|x\|, \|y\|\}\ge (8n+2)^{1/2}$.

To obtain a nontrivial estimate on $|L_n(x, y)|$ we follow the approach from \cite{PX3}.
Let $A_j^{(x)}:=-\frac{\partial}{\partial x_j}+x_j$.
Then the following identity holds (see the proof of Theorem~1 in \cite{PX3}):
\begin{align} \label{eq:hermite-1}
2^k(x_j-y_j)^k L_n(x,y) = \sum_{k/2 \le l \le k} c_{l,k}
\sum_{\nu=0}^{\infty} \Delta^l \ha\Big(\frac{\nu}{n}\Big)
\left(A_j^{(y)} - A_j^{(x)}\right)^{2l-k} \CH_\nu(x,y),
\end{align}
where $\Delta^l \ha(\frac{\nu}{n})$ is the $l$th forward difference applied to
the sequence $\{\ha(\frac{\nu}{n})\}_{\nu=0}^\infty$,
and the coefficients $c_{l,k}$ are given by
\begin{equation}\label{rep-c-lk}
    c_{l,k} = (-1)^{k-l} 4^{k-l} (2k-2l-1)!! \binom{k}{2l-k}
          = (-1)^{k-l} 2^{k-l} \frac{k!}{(k-l)!(2l-k)!}.
\end{equation}
Identity (\ref{eq:hermite-1}) follows from Lemma 3.2.3 in \cite{Th}; the constants $c_{l,k}$ are
given explicitly in \cite{PX3}.

Let us assume that $3\le k \le n/4$.
Using the estimate \cite[Lemma 3.2.2]{Th}
\begin{equation}\label{est-Hn}
\CH_n(x,x) \le c n^{d/2-1}, \qquad x \in \RR^d,
\end{equation}
one can follow the proof of Theorem 2.2 in \cite{PX3} to show that
\begin{align*}
& \Big| \left(A_j^{(y)} -  A_j^{(x)}\right)^{2l-k}   \CH_\nu(x,y)\Big|
   \le c (2\nu+4l-2k)^{(2l-k)/2}  \\
&  \qquad\qquad\qquad
     \times \sum_{i=0}^{2l-k} \binom{2l-k}{i}
            (\nu+i)^{(\frac{d}{2} -1)/2}    (\nu+2l-k+i)^{(\frac{d}{2} -1)/2} \\
&  \qquad\qquad
   \le 2^{(2l-k)/2} (\nu+2l-k)^{(2l-k +d-2)/2}
       \sum_{i=0}^{2l-k} \binom{2l-k}{i}      \\
& \qquad\qquad
   \le 2^{3 (2l-k)/2} (\nu+2l-k)^{(2l-k +d-2)/2} \le   2^{3 k/2} (\nu+k)^{(2l-k +d-2)/2}.
\end{align*}
Evidently,
$\left|\Delta^l \ha \left(\frac{\nu}{n}\right)\right|
\le n^{-l}\|\ha^{(l)}\|_{\infty}$
and  $\Delta^l \ha(\frac{\nu}{n})=0$ if $0\le \nu \le n/2 -l$ or $\nu\ge 2n$.
Using the above and the representation of $c_{l,k}$ from (\ref{rep-c-lk})
in \eqref{eq:hermite-1} we obtain
\begin{align*}
 2^k|x_j-y_j|^k | L_n(x,y)| & \le c 2^{3k/2} \sum_{k/2\le l \le k} \sum_{\nu = n/2-k}^{2n}
         (\nu+k)^{(2l-k +d-2)/2} |c_{l,k}| n^{-l} \|\ha^{(l)}\|_\infty\\
      & \le c 2^{2k} 3^{k+d} n^{(d-k)/2} \sum_{k/2\le l \le k}
               \frac{k!}{(k-l)! (2l-k)!} \|\ha^{(l)}\|_\infty,
\end{align*}
where we used that $\nu+k \le 2 n + n/4<3n$.
This along with estimate \eqref{est-derivative}
in Theorem \ref{thm:cutoff} leads to
\begin{align} \label{eq:hermite-2}
 2^k|x_j-y_j|^k   | L_n(x,y)| & \le c (12)^k n^{(d-k)/2} \sum_{k/2\le l \le k}
               \frac{k!}{(k-l)! (2l-k)!}  (\cc/\eps)^l l^l(\ln l)^{l(1+\eps)} \notag\\
    & \le c   n^{(d-k)/2} (12 \cc/\eps)^k (\ln k)^{k(1+\eps)} k!
              \sum_{k/2\le l \le k}  \frac{k^l}{(k-l)! (2l-k)!}.
\end{align}
We next establish the following estimate for the last sum above:
\begin{align} \label{eq:hermite-3}
      \sum_{k/2\le l \le k}  \frac{k^l}{(k-l)! (2l-k)!} \le 2 (2 e^2)^k.
\end{align}
To prove this we split the sum into two sums.
For $k/2\le l \le 3k/4$, by Stirling's formula, we have
$(k-l)! \ge [(k-l)/e]^{k-l} \ge  [ k/(4 e)]^{k-l}$ and hence
\begin{align*}
   \sum_{k/2\le l \le 3k/4}  \frac{k^l}{(k-l)! (2l-k)!} & \le
         (4 e)^{k/2}    \sum_{k/2\le l \le 3k/4}  \frac{k^{2l-k}}{(2l-k)!}  \\
  &  \le (4 e)^{k/2}    \sum_{m=0}^\infty   \frac{k^m}{m!}
    =  (4 e^3)^{k/2},
\end{align*}
whereas for $3k/4 < l \le k$, we have $(2l-k)! \ge  [k/(2e)]^{2l-k}$,
and hence
\begin{align*}
\sum_{3k/4< l \le k}  \frac{k^l}{(k-l)! (2l-k)!}
\le (2 e)^{k} \sum_{3k/4 < l \le k}\frac{k^{k-l}}{(k-l)!}\le  (2 e^{2})^k,
\end{align*}
and \eqref{eq:hermite-3} follows.
Substituting \eqref{eq:hermite-3} in \eqref{eq:hermite-2} we get
$$
2^k|x_j-y_j|^k   | L_n(x,y)| \le c n^{(d-k)/2} c_\eps^k  k^k(\ln k)^{k(1+\eps)},
\quad 1\le j \le k,
$$
where $c_\eps = 24 e^2 \cc/\eps\ge 24 e^2$; we used that $k! < k^k$.
Therefore,
\begin{equation} \label{eq:hermite-4}
|L_n(x,y)| \le c n^{d/2} \left( \frac{c_\eps k (\ln k)^{1+\eps}}{n^{1/2}\|x-y\|} \right)^k,
\quad 3 \le k \le n/4.
\end{equation}


We use (\ref{est-Hn}) and the Cauchy-Schwartz inequality to obtain the following trivial estimate
\begin{align}\label{est-Hermite-Ln}
|L_n(x,y)|
& \le c \sum_{j=0}^{2n}\Big(\sum_{|\a|=j}|\CH_\a(x)|^2\Big)^{1/2}
                       \Big(\sum_{|\a|=j}|\CH_\a(x)|^2\Big)^{1/2}\\
& = c \sum_{j=0}^{2n} \CH_j(x,x)^{1/2} \CH_j(y,y)^{1/2}\notag
\le cn^{d/2}.
\end{align}


Let $\max\{\|x\|, \|y\|\}\le (8n+2)^{1/2}$.
From (\ref{est-Hermite-Ln}) it follows immediately that estimate \eqref{Hermite-bound1} holds if
$
n^{1/2}\|x-y\|\le 3 c_\eps[\ln(e+n^{1/2}\|x-y\| )]^{1+\eps}.
$

Assume now that
$
n^{1/2}\|x-y\|> 3 c_\eps [\ln(e+n^{1/2}\|x-y\| )]^{1+\eps}
$
and choose
$$
k: = \left \lfloor \frac{n^{1/2}\|x-y\|}{c_\eps [\ln (e+ n^{1/2}\|x-y\|)]^{1+\eps} }\right
      \rfloor .
$$
Then $3\le k \le n/4$ as  $\|x-y\| \le 2(8n+2)^{1/2}$ and evidently
$
\frac{c_\eps k (\ln k)^{1+\eps}}{n^{1/2}\|x-y\|}\le e^{-1}.
$
Now, estimate \eqref{Hermite-bound1} follows
by (\ref{eq:hermite-4}) with $\CC = c_\eps^{-1}=\eps(24e^2\cc)^{-1}$.
\qed

\medskip

\noindent
{\bf Sub-exponentially localized needlets in the context of Hermite functions.}
For this construction, we shall utilize the cubature formula from \cite{PX3}
with nodes in $\cX_j\subset \RR^d$ consisting of $O(4^{jd})$ points on $\RR^d$
obtained as a product of nodal sets of univariate Gaussian quadrature formulas
and with positive coefficients $\{c_\xi\}_{\xi\in\cX_j}$, which is exact
for functions in $\bigoplus_{m=0}^N \CV_m$ with $N:=c4^{j+1}$
(for more details, see \cite{PX3}).
As before, we choose an admissible cutoff function $\ha$ of type (c)
obeying $(\ref{est-derivative})$ and define
$$
\psi_\xi(x):=c_\xi^{1/2}L_n(\xi, x),
\quad \xi\in\cX_j,
$$
where $L_n(x, y)$ with $n:=4^{j-1}$ is the kernel from (\ref{Hermire-kernels}).
Setting $\cX:=\cup_{j\ge 0}\cX_j$
we define the Hermite needlet system by
$\Psi:=\{\psi_\xi\}_{\xi\in\cX}$.
The sub-exponential localization of the Hermite needlets is inherited from
Theorem~\ref{thm:Hermite-localization-1}: For $\xi\in\cX_j$, $j\ge 0$,
\begin{equation}\label{local-Needlet-Hermite1}
|\psi_\xi(x)|
\le  c 2^{jd} \exp\Big\{-\frac{c'\eps 2^j \|x-\xi\|} {[\ln(e+2^j\|x-\xi\| )]^{1+\eps}}\Big\},
\qquad x\in  \RR^d,
\end{equation}
and
\begin{equation}\label{local-Needlet-Hermite2}
|\psi_\xi(x)|
\le ce^{-c''\|x\|^2}
\quad \mbox{if}\;\; \|x\|\ge c_*2^j
\end{equation}
for an appropriate constant $c_*$.
This is an improvement compared with the localization of the needlets from \cite{PX3}.

\section{Sub-exponentially localized kernels and frames induced by Laguerre functions}
\label{Laguerre-kernels}
\setcounter{equation}{0}

The Laguerre polynomials, defined by
$$
  L_n^\a(t) = \frac{1}{n!} t^{-\a} e^t \frac{d^n}{d t^n} (t^{n+\a} e^{-t}), \quad \a > -1,
  \quad n = 0,1,\ldots.
$$
are orthogonal on $\RR_+ =(0, \infty)$ with weight $t^\a e^{-t}$.
There are three types of Laguerre functions considered in the literature (see \cite{Th}),
defined by
\begin{equation}\label{def.Psi-n}
\CF_n^\alpha(t):=
 \Big(\frac{2\Gamma(n+1)}{\Gamma(n+\a+1)}\Big)^{1/2}e^{-t^2/2}L_n^\a(t^2),
\end{equation}
\begin{equation}\label{def.L-n}
\CL_n^\alpha(t):= \Big(\frac{\Gamma(n+1)}{\Gamma(n+\a+1)}\Big)^{1/2}
    e^{-t/2}t^{\a/2}L_n^\a(t),
\end{equation}
and
\begin{equation}\label{def.M-n}
\CM_n^\alpha(t):= (2t)^{1/2}\CL_n^\alpha(t^2).
\end{equation}
It is well known that $\{\CF_n^\alpha\}_{n\ge 0}$ is an orthonormal basis
for the weighed space $L^2(\RR_+, t^{2\a+1})$,
while $\{\CL_n^\alpha\}_{n\ge 0}$ and $\{\CM_n^\alpha\}_{n\ge 0}$
are orthogonal bases for $L^2(\RR_+)$.
Here we only consider the Laguerre functions $\{\CF_n^\alpha\}$ for $\a\ge 0$.
Analogous results for  $\{\CL_n^\alpha\}_{n\ge 0}$ and $\{\CM_n^\alpha\}_{n\ge 0}$
follow immediately as in \cite{KPPX}.

The $d$-dimensional tensor product Laguerre functions associated to
$\{\CF_n^\alpha\}$ are defined by
$$
\CF_\nu^\a (x): = \CF_{\nu_1}^{\a_1}(x_1) \cdots \CF_{\nu_d}^{\a_d}(x_d),
\quad x=(x_1, \dots , x_d),
$$
where $\nu = (\nu_1,\ldots,\nu_d) \in \NN_0^d$ and $\a = (\a_1,\ldots, \a_d)$.
The kernel of the orthogonal projector onto
$\CV_n:={\rm span} \{\cF^\a_\nu: |\nu|=n\}$ is given by
$$
\CF^\a_n(x, y):= \sum_{|\nu|=n} \CF^\a_\nu(x)\CF^\a_\nu(y).
$$

As elsewhere in this paper, we are interested in constructing sup-exponential localized
kernels of the form
\begin{equation}\label{Laguerre-Ln}
 L_n^\a(x,y):= \sum_{j=0}^\infty \ha\Big(\frac{j}{n}\Big)
    \CF_j^\a(x,y), \quad x, y \in \RR_+^d,
\end{equation}
where $\ha$ is an admissible cutoff function (see Definition~\ref{cutoff-d1}).
Note that in \cite{KPPX} it is proved that for admissible cutoff functions $\ha$ the kernels
$L_n^\a(x,y)$ have faster than the reciprocal of any polynomial decay away from
the main diagonal in $\RR_+^d\times\RR_+^d$.

We shall use some of the notation and results from \cite{KPPX}.
Recall our standing notation for norms in $\RR^d$:
$\|x\|:=\max_i |x_i|$, $\|x\|_2:= (\sum_i |x_i|^2)^{1/2}$, and
$|x|=\|x\|_1:=\sum_i |x_i|$.
We shall also need the quantity
$$
W_\a(n;x):= \prod_{j=1}^d (x_j+n^{-1/2})^{2\alpha_j+1}, \qquad x \in  \RR^d_+.
$$


\begin{theorem}\label{thm:Laguerre-localization-1}
Let $0<\eps \le 1$ and assume that  $\ha$ in $(\ref{Laguerre-Ln})$
is an admissible cutoff function obeying inequality $(\ref{est-derivative})$
in Theorem~\ref{thm:cutoff}.
Then for $x,y \in \RR_+^d$,
\begin{equation}\label{Laguerre-bound1}
|L_n^\a (x,y)|
\le \frac{c n^{d/2}}{\sqrt{W_\a(n; x)}\sqrt{W_\a(n; y)}}
\exp\Big\{-\frac{c'\eps n^{1/2}\|x-y\|}{[\ln(e+n^{1/2}\|x- y\|)]^{1+\eps}}\Big\},
\end{equation}
where
$c'>0$ is a constant depending only on $d$ and $\alpha$,
and $c$ depends on $d$, $\alpha$, $\eps$.

As in similar situations before, for an appropriate cutoff function $\ha$ the term
$[\ln(e+n^{1/2}\|x- y\|)]^{1+\eps}$ above can be replaced by any product of multiple logarithms
as in $(\ref{product2})$.

In addition, for some constant $c''>0$
\begin{equation}\label{Laguerre-bound2}
|L_n^\a (x,y)|
\le ce^{-c''\max\{\|x\|^2, \|y\|^2\}}
\quad \mbox{if}\; \max\{\|x\|, \|y\|\}\ge (12n+3\|\alpha\|+3)^{1/2}.
\end{equation}
\end{theorem}

Besides Theorem~\ref{thm:cutoff} a main new ingredient in the proof of this theorem
will be the following result for Laguerre polynomials:


\begin{proposition}\label{prop:Laguerre-bound}
For $n\ge 1$ and $-1/2\le\a \le n$, we have
\begin{equation} \label{Laguerre-est*}
| L_n^\a(t)| e^{-t/2} \le c2^\a (n/t)^{\alpha/2}, \quad 0<t<\infty,
\end{equation}
where $c>0$ is an absolute constant.
The lower bound $-1/2$ for $\a$ can be replaced by any constant $\alpha_0>-1$
but $c$ will depend on $\alpha_0$.
\end{proposition}

\noindent
{\bf Proof.}
We shall utilize some estimates on $|L_n^\a(t)|$ established in \cite{Kra}.
As is shown in \cite[Theorem 2]{Kra} for $\a\ge 24$ and $n\ge 1$
\begin{align}\label{Kra1}
\max_{0 < t \le \frac{2(\a+1)^2}{2n+\a+1}} t^{\a+1} e^{-t} [ L_n^\a(t) ]^2
 < 650 \frac{\Gamma(n+\alpha+1)}{\Gamma(n+1)}  \frac{(\a+1)^{4/3}}{n^{1/6}(n+\a+1)^{5/6}}.
\end{align}
For the next estimate, let us denote briefly
$$
s_n := \sqrt{n+\a+1} + \sqrt{n} \quad \hbox{and} \quad
q_n := \sqrt{n+\a+1} - \sqrt{n}.
$$
Then, for $\a \ge 24$ and $n \ge 35$
\cite[Theorem  4]{Kra},
\begin{align}\label{Kra3}
t^{\a+1} e^{-t} [ L_n^\a(t) ]^2
< 680 \frac{\Gamma(n+\alpha+1)}{\Gamma(n+1)} \frac{t}{\sqrt{ (t-q_n^2)(s_n^2 - t)}}
\quad  \hbox{for  $t \in (q_n^2, s_n^2)$}.
\end{align}

We shall consider several cases for $\a$, $n$, and $t$.

\smallskip


{\em Case} 1: $-1/2\le \alpha<24$ and $n\ge 1$ or $1\le n < 35$ (then $\alpha \le n <35$).
In this case, estimate (\ref{Laguerre-est*}) ($0<t<\infty$)
follows from the standard estimate in \cite[\S 8.22]{Sz} with a constant depending on $\alpha$
(see also \cite[Lemma 1.5.3]{Th} or \cite[(2.11)]{PX3}).

\smallskip


{\em Case} 2: $n\ge 35$, $24\le \alpha\le n$, and $0 < t \le 1/n$.
A standard upper bound for $|L_n^\a(t)|$ ($\a\ge 0$) is given by (see e.g. \cite[(3)]{Kra})
\begin{equation} \label{triv-Laguerre}
|L_n^\a(t)|e^{-t/2}\le \frac{(\a+1)_n}{n!}
= \frac{\Gamma(n+\alpha+1)}{\Gamma(\a+1)\Gamma(n+1)},
\quad 0<t<\infty.
\end{equation}
Note that by Stirling's formula
$\Gamma(t+1) \sim (t/e)^t\sqrt{t}$ on $[1, \infty)$ and hence
\begin{equation} \label{Gamma-est}
\frac{\Gamma(n+\alpha+1)}{\Gamma(n+1)}
\le c\frac{((n+\a)/e)^{n+\a}\sqrt{n+\a}}{(n/e)^n\sqrt{n}}
\le c2^\a n^\a,
\end{equation}
where $c>0$ is an absolute constant.
Using this in (\ref{triv-Laguerre}) we get
$$
|L_n^\a(t)|e^{-t/2} \le cn^\a\le c(n/t)^{\a/2},
\quad 0<t\le 1/n,
$$
which gives (\ref{Laguerre-est*}) in the case under consideration.

\smallskip


{\em Case} 3: $n\ge 35$, $24\le \alpha\le n$, and $\frac{1}{n} < t \le \frac{2(\a+1)^2}{2n+\a+1}$.
From (\ref{Kra1}) and (\ref{Gamma-est}) we get
$$
e^{-t}[ L_n^\a(t) ]^2
\le\frac{c(\a+1)^{4/3}}{n^{1/6}(n+\a+1)^{5/6}}\frac{2^\a n^\a}{t^{\a+1}}
\le \frac{c\alpha^{4/3}2^\a}{nt} \Big(\frac{n}{t}\Big)^{\a}
\le c\alpha^{4/3}2^\a \Big(\frac{n}{t}\Big)^{\a},
$$
which implies (\ref{Laguerre-est*}).

\smallskip


{\em Case} 4: $\frac{2(\a+1)^2}{2n+\a+1}<t\le n$.
We have
$q_n^2 < \frac{(\a+1)^2}{4n}< \frac{1}{2}\frac{2(\alpha+1)^\a}{2n+\a+1}$
and $s_n^2> 4n$.
We use these and (\ref{Gamma-est}) in (\ref{Kra3}) to obtain
$$
t^{\a} e^{-t} [ L_n^\a(t) ]^2
\le \frac{c2^\a n^\a}{\sqrt{ (t-q_n^2)(s_n^2 - t)}}
\le \frac{c2^\a n^\a}{\sqrt{n/(n+\a+1)}}
\le c2^\a n^\a,
$$
which again gives (\ref{Laguerre-est*}).
The proof of the proposition is complete.
\qed

\medskip

\noindent
{\bf Proof of Theorem~\ref{thm:Laguerre-localization-1}.}
Note first that from Lemma 1.5.3 in \cite{Th} or \cite[\S8.22]{Sz}
it follows that
\begin{equation} \label{Est-Lag}
|L_n^\a(t)| e^{-t/2} \le  cn^\a
\quad\mbox{for} \; t\ge 0
\quad\mbox{and}\quad
|L_n^\a(t)| e^{-t/2} \le  c e^{-\g n}
\quad\mbox{for} \; t \ge 3N/2,
\end{equation}
where $\gamma>0$ is a constant and $N:=4n+2\alpha+2$.
These immediately lead to the estimates (see also \cite{KPPX}):
\begin{equation} \label{Est-Fn}
\|\cF_n^\a\|_\infty\le cn^{\a/2}
\quad\mbox{and}\quad
|\cF_n^\a(t)| \le ce^{-\gamma't^2}
\quad\mbox{for} \; t \ge (3N/2)^{1/2},
\end{equation}
which readily imply (\ref{Laguerre-bound2}).
In turn estimate (\ref{Laguerre-bound2}) obviously implies that (\ref{Laguerre-bound1})
holds if
$\max\{\|x\|, \|y\|\}\ge (12n+3\|\alpha\|+3)^{1/2}$.

The proof of (\ref{Laguerre-bound1}), when $\max\{\|x\|, \|y\|\}< (12n+3\|\alpha\|+3)^{1/2}$,
will rely on the developments in \cite{KPPX},
where in particular it is shown that $L_n^\a(x,y)$
has the representation (see \cite[(8.3)-(8.4)]{KPPX})
\begin{align}\label{Laguerre-1}
L_n^\a(x, y)
&= c  \int_{[0,\pi]^d } \CK_n^\a
     \Big (\|x\|_2^2+\|y\|_2^2+2\sum_{i=1}^d x_iy_i\cos\theta_i \Big)\\
& \qquad \times \prod_{i=1}^d j_{\a_i-1/2}(x_iy_i\cos\theta_i)
      \sin^{2\a_i}\theta_i \,d\theta, \notag
\end{align}
where $j_\a(t):=t^{-\a}J_\a(t)$ with $J_\a(t)$ being the Bessel function,
and the kernel $\CK_n^\a$ is given by
\begin{align}\label{Laguerre-2}
\CK_n^\a(t)
 =  \sum_{m=0}^\infty \Delta^{k+1}\ha\Big(\frac{m}{n}\Big)
L_m^{|\a|+k+d} (t) e^{- t/2}. 
\end{align}
Here $k\ge 0$ is arbitrary and the finite difference
$\Delta^{k+1}$ is with respect to $m$.

We~claim that $\CK_n^\a$ satisfies the estimate
\begin{equation} \label{Laguerre-3}
| \CK_n^\a(t)| \le c n^{|\a| + d}
\exp\left\{-\frac{\CC(n t)^{1/2}}{[\ln(e+(nt)^{1/2})]^{1+\eps}}\right\},
\quad t \le 4d(12n+3\|\alpha\|+3),
\end{equation}
where  $\CC=c'\eps$ with $c'>0$ independent of $\eps$.

Suppose $t \le 4d(12n+3\alpha+3)$ and let $2 \le k \le n/4$.
Evidently, we have
$
\left|\Delta^{k+1}\ha \left(\frac{\nu}{n}\right)\right|
\le n^{-k-1}\|\ha^{(k+1)}\|_{\infty}
$
and
$\Delta^{k+1} \ha(\frac{\nu}{n})=0$ if $0\le \nu \le n/2 -k-1$ or $\nu\ge 2n$.
From these, estimate \eqref{est-derivative} from Theorem \ref{thm:cutoff}, and
(\ref{Laguerre-est*}) it follows that
\begin{align} \label{Laguerre-4}
|\CK^\a_n(t)|
& \le  c\sum_{n/2 - k -1< m < 2n}
\frac{1}{n^{k+1}}\Big(\frac{4m}{t}\Big)^{(|\a| +k+d)/2} \|\ha^{(k+1)}\|_\infty \notag\\
& \le cn^{-k} \Big(\frac{8n}{t}\Big)^{(|\a| +k+d)/2}
        (\cc/\eps)^{k+1} (k+1)^{k+1} [\ln (k+1)]^{(k+1)(1+\eps)} \\
& \le c  n^{|\a|+d} \left(\frac{8}{nt}\right)^{(|\a| +d-1)/2}
         \left[ \frac{c_\eps (k+1) [\ln (k+1)]^{1+\eps}} {(nt)^{1/2}} \right]^{k+1}, \notag
\end{align}
where $c_\eps = 64\cc/\eps$, $\cc\ge 1$.


Exactly as above, but using the first estimate in (\ref{Est-Lag}) instead of (\ref{Laguerre-est*})
and (\ref{Laguerre-2}) with $k=0$ we get
\begin{equation}\label{Laguerre-5}
|\CK^\a_n(t)|
\le c \sum_{n/2 -1< m < 2n}n^{-1}m^{|\a| +d} \|\ha'\|_\infty
\le cn^{|\a|+d}.
\end{equation}

Suppose
$(nt)^{1/2} \ge 2c^\flat c_\eps [\ln(e+\sqrt{nt})]^{1+\eps}$,
where $c^\flat:= (12d(5+\|\a\|))^{1/2}$, and choose
$$
k: = \left\lfloor \frac{\sqrt{nt}} {c^\flat c_\eps [\ln(e+\sqrt{nt})]^{1+\eps}} \right \rfloor -1.
$$
Then $3\le k+1 \le n/4$ as $t \le 4d(12n+3\|\alpha\|+3)\le (c^\flat)^2n$,
and it is easy to see that
$\frac{c_\eps (k+1) [\ln (k+1)]^{1+\eps}} {\sqrt{nt}} \le e^{-1}$.
Hence, by (\ref{Laguerre-4}) it follows that (\ref{Laguerre-3}) holds
with $\CC = (c^\flat c_\eps)^{-1}$,
where we used that $nt\ge (c^\flat c_\eps)^2$.

If $(nt)^{1/2} < c^\flat c_\eps [\ln(e+\sqrt{nt})]^{1+\eps}$, then \eqref{Laguerre-3}
is immediate from (\ref{Laguerre-5}).
This completes the proof of (\ref{Laguerre-3}).


To prove \eqref{Laguerre-bound1}, we will use \eqref{Laguerre-3}
with  $t = \|x\|_2^2+\|y\|_2^2+2\sum_{i=1}^d x_iy_i\cos\theta_i$.
Evidently,  $t \ge \|x - y\|^2$ and
$t\le d(\|x\|+\|y\|)^2\le 4d(12n+3\|\a\|+3)$ as in (\ref{Laguerre-3}).
Also,
$
\exp\big\{\frac{-\frac{1}{2}c^\diamond x}{[\ln(e+x)]^{1+\eps}}\big\}
\le c(1+x^2)^{-2|\a|-d},
$
$x>0$, where $c>0$ is a sufficiently large constant.
Therefore, from \eqref{Laguerre-3} it follows that
$$
| \CK_n^\a(t) | \le  c  n^{|\a|+d}
\exp\left\{-\frac{\frac{1}{2}\CC n^{1/2}\|x-y\|}{[\ln(e+n^{1/2}\|x-y\|)]^{1+\eps}}\right\}
\left( \frac{1}{1+ nt} \right)^{\tau}
$$
with $\tau = 2|\a| + d$. Combining this with \eqref{Laguerre-1}
and the fact that $|j_{\a-1/2}(t)| \le c_\a$ for $t \in \RR_+$,
$\a\ge 0$, (see \cite{KPPX}) we arrive at the estimate
\begin{align*}
|L_n^\a(x, y)| & \le c  n^{|\a|+d}
      \exp\left\{-\frac{\frac{1}{2}\CC n^{-1/2}\|x-y\|}{[\ln(e+n^{1/2}\|x-y\|)]^{1+\eps}}\right\} \\
&   \qquad \times  \int_{[0,\pi]^d}  \frac{ \prod_{i=1}^d
      \sin^{2\a_i}\theta_i \,d\theta}
{\left[1+ n\Big (\|x\|_2^2+\|y\|_2^2+2\sum_{i=1}^d x_iy_i\cos\theta_i \Big) \right]^\tau}.\notag
\end{align*}
Denote the last integral above by $J$.
We now estimate $J$ quite similarly as in the proof of Theorem~3.2 in \cite{KPPX}.
Substituting $\theta_i:=\pi-t_i$ and using that
$1-\cos t=2\sin^2\frac{t}{2}\sim t^2$ we get
\begin{align*}
J&\le c\int_{[0,\pi]^d}  \frac{ \prod_{i=1}^d t^{2\a_i} \,dt}
{\left[1+ n\Big (\|x-y\|^2+\sum_{i=1}^d x_iy_it_i^2 \Big) \right]^\tau}\\
&\le c\int_{[0,\pi]^\ell}  \frac{ \prod_{i=1}^\ell t^{2\a_i} \,dt}
{\left[1+ n\Big (\|x-y\|^2+\sum_{i=1}^\ell x_iy_it_i^2 \Big) \right]^\tau},
\quad 0\le \ell\le d.
\end{align*}
For a fixed $1\le \ell\le d$ denote $\alpha':=(\a_1, \dots, \a_\ell)$.
Then substituting $u_i=t_i(x_iy_i)^{1/2}$ above, we get
\begin{align*}
J & \le \frac{c}{n^{|\a'|}\prod_{i=1}^\ell (x_i y_i)^{\a_i+1/2}}  \prod_{i=1}^\ell
\int_0^{\pi(x_iy_i)^{1/2}}\frac{du}{\big[1+n(\|x-y\|^2+\sum_{i=1}^\ell u_i^2)\big]^{\tau-|\a'|}},
\end{align*}
where we used $|\a'|$ terms from the denominator to cancel the numerator.
We now enlarge the domain of integration to $\RR^\ell$ and use spherical coordinates
to bound the product of integrals above by
\begin{align*}
\int_0^\infty \frac{r^{\ell-1} dr } {\big[1+n(\|x-y\|^2+r^2)\big]^{\tau-|\a'|}}
\le  \frac{c} {n^{\ell/2}\big(1+n\|x-y\|^2\big)^{\tau-|\a'|-\ell/2}},
\end{align*}
where we used that $\tau:=2|\a|+d>|\a|+d/2$.
Therefore, we have for $0\le \ell\le d$ (the case $\ell=0$, $|\a'|:=0$, is trivial)
\begin{align*}
n^{|\a|+d/2}J
& \le \frac{cn^{|\a|+d/2}}{n^{|\a'|+\ell/2}
\prod_{i=1}^\ell (x_i y_i)^{\a_i+1/2}(1+n\|x-y\|^2\big)^{\tau-|\a|-d/2}}\\
&\le \frac{c}{\prod_{i=1}^\ell (x_i y_i)^{\a_i+1/2}
              \prod_{i=\ell+1}^d (n^{-1})^{\a_i+1/2}
(1+n^{1/2}\|x-y\|\big)^{|\a|+d/2}},
\end{align*}
which readily implies
\begin{align*}
n^{|\a|+d/2}J
\le \frac{c}{\prod_{i=1}^\ell (x_i y_i+n^{-1})^{\a_i+1/2}
(1+n^{1/2}\|x-y\|\big)^{|\a|+d/2}}.
\end{align*}
Using now the simple inequalities
$$
(x_i+n^{-1/2})(y_i+n^{-1/2}) \le 3(x_iy_i+n^{-1})(1+n^{1/2}\|x-y\|), \quad 1\le i\le d,
$$
we obtain
$$
n^{|\a|+d/2}J\le \frac{c}{\sqrt{W_\a(n; x)}\sqrt{W_\a(n; y)}}
$$
and estimate \eqref{Laguerre-bound1} follows.
\qed

\medskip

\noindent
{\bf Sub-exponentially localized needlets induced by Laguerre functions.}
The construction here is similar as in previous sections, in particular,
see the case of Hermite functions.
We start with the cubature formula from \cite{KPPX}
with nodes in $\cX_j\subset \RR_+^d$ consisting of $O(4^{jd})$ points in $\RR^d$
obtained as product of the nodal sets of univariate Gaussian quadrature formulas
on $\RR_+$ and with positive coefficients $\{c_\xi\}_{\xi\in\cX_j}$, which is exact
for functions in $\bigoplus_{m=0}^N \CV_m$ with $N:=c4^{j+1}$ (see \cite{KPPX}).
We next choose an admissible cutoff function $\ha$ of type (c)
obeying (\ref{est-derivative}) and then define
$$
\psi_\xi(x):=c_\xi^{1/2}L_n(\xi, x),
\quad \xi\in\cX_j,
$$
where $L_n(x, y)$ with $n:=4^{j-1}$ is the kernel from (\ref{Hermire-kernels}).
Notice that here the dilation on the frequency side is by a factor of 4.
Setting $\cX:=\cup_{j\ge 0}\cX_j$,
we define the Laguerre needlet system by
$\Psi:=\{\psi_\xi\}_{\xi\in\cX}$.
The sub-exponential localization of the Laguerre needlets
follows by Theorem~\ref{thm:Laguerre-localization-1} and takes the form:
For $\xi\in\cX_j$, $j\ge 0$,
\begin{equation}\label{local-Needlet-Laguerre1}
|\psi_\xi(x)|
\le \frac{c 2^{jd}}{\sqrt{W_\a(n; \xi)}}
\exp\Big\{-\frac{c'\eps 2^j\|x-\xi\|}{[\ln(e+2^j\|x-\xi\|)]^{1+\eps}}\Big\},
\quad x\in\RR_+^d,
\end{equation}
and
\begin{equation}\label{local-Needlet-Laguerre2}
|\psi_\xi(x)|
\le ce^{-c''\|x\|^2}
\quad \mbox{if}\; \|x\|\ge c_*2^j
\end{equation}
for an appropriate constant $c_*$.
This is a natural improvement of the localization of the needlets from \cite{KPPX}.

\section{The localization principle from  \S\ref{Examples} fails for product Jacobi polynomials}
\label{not-valid-principle}
\setcounter{equation}{0}

In this section, we show that in contrast to the cases considered in previous sections
surprisingly the localization principle described in \S\ref{Examples}
is no longer valid for 2-d tensor product Legendre or Chebyshev polynomials
and products of Legendre and Chebyshev polynomials.

The $n$th degree Legendre polynomial $P_n$ is defined by
$$
P_n(x):= \frac{1}{n!2^n} \Big(\frac{d}{dx}\Big)^n(x^2-1)^n
$$
and it is known that $P_n(1)=1$, $P_n(-1)=(-1)^n$
and $\int_{-1}^1P_n^2(x)dx=(n+1/2)^{-1}$.
Denote $\tP_n:=(n+1/2)^{1/2}P_n$ and note that
$\{\tP_n\}_{n\ge 0}$ is an orthonormal basis for $L^2[-1, 1]$.

The 2-d tensor product Legendre polynomials are defined by
$$
\tP_\nu(x):= \tP_{\nu_1}(x_1)\tP_{\nu_2}(x_2),
\quad \nu=(\nu_1, \nu_2), \; x=(x_1, x_2).
$$
For an admissible univariate cutoff function $\ha\ge 0$ (see Definition~\ref{cutoff-d1})
define
\begin{equation} \label{Legendre-kernel}
L_n(x, y):= \sum_{m=0}^\infty \ha\Big(\frac{m}{n}\Big) \tP_m(x, y),
\quad\mbox{where}\quad
\tP_m(x, y):= \sum_{|\nu|=m}\tP_\nu(x)\tP_\nu(y).
\end{equation}
Here $|\nu|:=\nu_1+\nu_2$ and $x, y\in [-1, 1]^2$.
In more detail
$$
\tP_m(x, y)
= \sum_{\nu_1+\nu_2=m}
(\nu_1+1/2)(\nu_2+1/2)P_{\nu_1}(x_1)P_{\nu_2}(x_2)P_{\nu_1}(y_1)P_{\nu_2}(y_2).
$$
In order to show that $L_n(x, y)$ has very poor localization at some points
on $[-1, 1]^2\times[-1, 1]^2$, which are away from the main diagonal, we fix $y=(1,1)$.
Using that $P_n(-1)=(-1)^n$ one easily verifies the identity
\begin{equation} \label{Legendre-eq1}
\tP_m(1, -1, 1, 1) = \sum_{j=0}^m (-1)^j(j+1/2)(m-j+1/2)=(1+(-1)^m)/8.
\end{equation}
Indeed, for odd $m$ this follows from the symmetry of the terms in the sum
and in the case $m=2k$ it follows by simple manipulations from the case $m=2k-1$ .

By (\ref{Legendre-kernel})-(\ref{Legendre-eq1}) it follows that
$$
L_n(1, -1, 1, 1) = \frac{1}{4}\sum_{j=0}^\infty\ha\Big(\frac{2j}{n}\Big)
= \frac{n}{8}\int_{0}^\infty \ha(t)dt+\frac{\ha(0)}{8}+\cO(n^{-1}),
$$
which shows that $L_n(x, y)$ has no localization whatsoever for $x=(1,-1)$, $y=(1,1)$.
The behavior of the sequence $\{L_n(x, (1,1))\}_{n=0}^\infty$ is similar
when $x$ is any other point on the lines $\{(x_1,1)\}$ and $\{(1,x_2)\}$.
However, the kernels $L_n(x, y)$ are very well localized in large portion
of $[-1, 1]^2\times[-1, 1]^2$.

It is slightly more complicated to show that the situation is quite the same
for Chebyshev polynomials. In this case (\ref{Legendre-kernel}) holds with
$$
\tP_m(x, y)
= \frac{4}{\pi^2}\sum_{j=0}^{m}
\Big(1-\frac{\delta_{j,0}}{2}\Big)\Big(1-\frac{\delta_{j,m}}{2}\Big)
T_j(x_1)T_{m-j}(x_2)T_j(y_1)T_{m-j}(y_2),
$$
where $T_n(v)=\cos n\arccos v$ and $\delta_{j,k}$ is the Kroneker delta
(see the definition of Chebyshev polynomials given in \S\ref{Examples}).
Using that $T_n(-1)=(-1)^n$ one easily verifies the identity
\begin{equation*}
\tP_m(1, -1, 1, 1) = \frac{4}{\pi^2}\sum_{j=0}^m
\Big(1-\frac{\delta_{j,0}}{2}\Big)\Big(1-\frac{\delta_{j,m}}{2}\Big)(-1)^{m-j}
=\frac{\delta_{m,0}}{\pi^2}.
\end{equation*}
Hence
$$
L_n(1, -1, 1, 1) = \frac{\ha(0)}{\pi^2},
$$
which implies that there is no localization for cutoff functions of type (a)
in Definition~\ref{cutoff-d1}.

The above identity is inconclusive  for cutoff functions $\ha$ of type (b) and (c).
In fact, for every such $\ha$ and every segment
$$\{x(t)=(1-2t\cos\theta,1-2t\sin\theta)~:~t\in [0,\min\{\sec\theta,\csc\theta\}]\}$$
with fixed $\theta\in[0,\pi/2]$ the sequence $\{L_n(x(t), (1,1))\}_{n=0}^\infty$
has excellent localization when $t$ increases.
This localization, however, is not uniform over $\theta\in[0,\pi/2]$ and even
$\{L_n(x, (1,1))\}_{n=0}^\infty$ is not localized at all for $x\in[-1,1]^2$!
In order to demonstrate this fact we consider the function $F_n(x_1)=L_n(x_1,-1,1,1)$.
Straightforward calculations give
$$
F_n(x_1)
=\frac{\ha(0)}{\pi^2}-\sum_{j=1}^\infty
\frac{2}{\pi^2}\ha\Big(\frac{j}{n}\Big)\sin (j\arccos x_1) \tan \Big(\frac{1}{2}\arccos x_1\Big).
$$
In particular, $F_n$ is an algebraic polynomial of degree at most $\ell n$,
where $\supp \ha \subset [0,\ell]$. The above identity implies for its derivative
$$
F_n'(1) = \frac{4}{\pi^2}\sum_{j=1}^\infty \frac{j}{2}\ha\Big(\frac{j}{n}\Big)
= \frac{2n^2}{\pi^2}\int_{0}^\infty t\ha(t)dt+\cO(1).
$$
Now, because of $\int_{0}^\infty t\ha(t)dt>0$, Markov's inequality implies that
$$\|F_n\|_{L^\infty[-1,1]}>const,$$
which validates our assertion.

Another example is the weight $(1-x_1^2)^{-1/2}$ in $[-1,1]^2$,
which is associated with products of Chebyshev and Legendre polynomials.
In this case (\ref{Legendre-kernel}) holds with
$$
\tP_m(x, y)
= \frac{2}{\pi}\sum_{j=0}^{m}
\Big(1-\frac{\delta_{j,0}}{2}\Big)\Big(m-j+\frac{1}{2}\Big)T_j(x_1)T_{m-j}(x_2)P_j(y_1)P_{m-j}(y_2).
$$
Now
\begin{equation*}
\tP_m(1, -1, 1, 1) = \frac{2}{\pi}\sum_{j=0}^{m}
\Big(1-\frac{\delta_{j,0}}{2}\Big)\Big(m-j+\frac{1}{2}\Big)(-1)^{m-j}=\frac{(-1)^m}{2\pi}
\end{equation*}
and hence
$$
L_n(1, -1, 1, 1) = \frac{\ha(0)}{4\pi}+\cO(n^{-1}).
$$
As for 
tensor products of Chebyshev polynomials the localization of
the cutoff functions of type (a) is ruled out by this argument,
but for those of types (b) and (c) we consider the derivative of
the function $F_n(x_1)=L_n(x_1,-1,1,1)$.
Here we have
$$
F_n'(1) = \frac{2}{\pi}\sum_{j=1}^\infty \frac{j}{2}\Big(\ha\Big(\frac{2j-1}{n}\Big)
+\ha\Big(\frac{2j}{n}\Big)\Big)
= \frac{n^2}{2\pi}\int_{0}^\infty t\ha(t)dt+ \frac{n}{4\pi}\int_{0}^\infty \ha(t)dt+\cO(1)
$$
and hence
$$\|F_n\|_{L^\infty[-1,1]}>const.$$

The above facts lead us to the  {\bf conclusion} that
one cannot expect good localization of kernels of
the form (\ref{Legendre-kernel}) for tensor product Jacobi polynomials
or cross product bases.
This kind of bases apparently have completely different nature compared to
e.g. orthogonal polynomials on the simplex (\S\ref{simplex-kernels})
and the ball (\S\ref{ball-kernels}),
or multivariate Hermite and Laguerre functions
(\S\ref{Hermite-kernels}-\ref{Laguerre-kernels})
for which the localization principle is valid.

In the case of tensor product Jacobi polynomials truly multivariate cutoff
functions need to be employed.
This sort of cutoff functions and the associated kernels and needlets
will be developed in a follow up paper.

\end{document}